\documentclass[11pt,a4paper]{article}

\usepackage{epic,eepic,epsfig,amsmath,amssymb,amsthm}
\usepackage{t1enc}
\usepackage{color}

\usepackage{bbm}
\def\virgp{\raise 2pt\hbox{,}}

\renewcommand{\geq}{\geqslant}
\renewcommand{\leq}{\leqslant}
\def\N{{\mathbb N}}

\def\R{{\mathbb R}}

\def\virgp{\raise 2pt\hbox{,}}
\def\cdotpv{\raise 2pt\hbox{;}}

\def\1{\mathbbm{1}}

\newtheorem{theorem}{Theorem}[section]

\newtheorem{proposition}[theorem]{Proposition}

\newtheorem{pte}[theorem]{Property}

\theoremstyle{remark}
\newtheorem{remark}{Remark}[section]

\theoremstyle{definition}
\newtheorem{definition}{Definition}[section]

\newtheorem*{notation}{Notation}

\theoremstyle{definition}

\theoremstyle{definition}

\begin{document}

\title{\textbf{Resistance metric, and spectral asymptotics, on the graph of the Weierstrass function}}

\author{\LARGE{\textbf{Claire David}}}

\maketitle
\centerline{Sorbonne Universit\'es, UPMC Univ Paris 06}

\centerline{CNRS, UMR 7598, Laboratoire Jacques-Louis Lions, 4, place Jussieu 75005, Paris, France}

%\begin{abstract}

%\end{abstract}

\maketitle

\section{Introduction}

\hskip 0.5cm Following our work on the graph of the Weierstrass function~\cite{ClaireGB}, in the spirit of those of~J. Kigami ~\cite{Kigami1989},~\cite{Kigami1993}, and~\cite{Strichartz1999}, \cite{StrichartzLivre2006}, which enabled us to build a Laplacian on the aforementioned graph, it was natural to go further and give the related explicit resistance metric. In doing so, we made calculations that directly enable one to obtain the box dimension of the graph, in a simpler way than~\cite{Kaplan1984}or~\cite{HuLau1993}. \\

The aim of this work is twofold. We had a special interest in the study of the spectral properties of the Laplacian. In~\cite{ClaireGB}, we have given the explicit the spectrum on the graph of the Weierstrass function. In the case of Laplacians on post-critically finite fractals, previous works, by J.~Kigami and M.~Lapidus~\cite{KigamiLapidus1993}, and R.~S.~Strichartz~\cite{StrichartzLivre2006}, make the link between resistance metric, and asymptotic properties of the spectrum of the Laplacian, by means of an analoguous of Weyl's formula.\\

So we asked ourselves wether those results were still valid, for the graph of the Weierstrass function.\\

\section{Framework of the study}

\noindent In this section, we recall results that are developed in~\cite{ClaireGB}.

\vskip 1cm

\begin{notation}
In the following,~$\lambda$ and~$N_b$ are two real numbers such that:

$$0 <\lambda<1 \quad, \quad N_b\,\in\,\N \quad  \text{and} \quad \lambda\,N_b > 1 $$

\noindent We will consider the~($1-$periodic) Weierstrass function~${\cal W}$, defined, for any real number~$x$, by:

$$   {\cal W}(x)=\displaystyle \sum_{n=0}^{+\infty} \lambda^n\,\cos \left ( 2\,  \pi\,N_b^n\,x \right)
$$

\end{notation}
\vskip 1cm

We place ourselves, in the sequel, in the Euclidean plane of dimension~2, referred to a direct orthonormal frame. The usual Cartesian coordinates are~$(x,y)$.\\

\vskip 1cm

\noindent The restriction~$\Gamma_{\cal W}$ to~$[0,1[ \times \R$, of the graph of the Weierstrass function, is approximated by means of a sequence of graphs, built through an iterative process. To this purpose, we introduce the iterated function system of the family of~$C^\infty$ contractions from~$\R^2$ to~$\R^2$:
$$\left \lbrace T_{0},...,T_{N_b-1} \right \rbrace$$

\noindent where, for any integer~$i$ belonging to~$\left \lbrace 0,...,N_b-1  \right \rbrace$, and any~$(x,y)$ of~$\R^2$:
$$ T_i(x,y) =\left( \displaystyle \frac{x+i}{N_b}, \lambda\, y + \cos\left(2\,\pi \,\left(\frac{x+i}{N_b}\right)\right) \right)$$

\vskip 1cm

\begin{pte}

$$   \Gamma_{\cal W} =  \underset{  i=0}{\overset{N_b-1}{\bigcup}}\,T_{i}(\Gamma_{\cal W})$$

\end{pte}

\vskip 1cm

\begin{definition}
\noindent For any integer~$i$ belonging to~$\left \lbrace 0,...,N_b-1\right \rbrace  $, let us denote by:

$$P_i=(x_i,y_i)=\left(\displaystyle \frac{i}{N_b-1},\displaystyle\frac{1}{1-\lambda}\,\cos\left ( \displaystyle\frac{2\,\pi\,i}{N_b-1}\right ) \right) $$

\noindent the fixed point of the contraction~$T_i$.\\

\noindent We will denote by~$V_0$ the ordered set (according to increasing abscissa), of the points:

$$\left \lbrace P_{0},...,P_{N_b-1}\right \rbrace$$

\noindent The set of points~$V_0$, where, for any~$i$ of~\mbox{$\left \lbrace  0,...,N_b-2  \right \rbrace$}, the point~$P_i$ is linked to the point~$P_{i+1}$, constitutes an oriented graph (according to increasing abscissa)), that we will denote by~$ \Gamma_{{\cal W}_0}$.~$V_0$ is called the set of vertices of the graph~$ \Gamma_{{\cal W}_0}$.\\

\noindent For any natural integer~$m$, we set:
$$V_m =\underset{  i=0}{\overset{N_b-1}{\bigcup}}\, T_i \left (V_{m-1}\right )$$

\noindent The set of points~$V_m$, where two consecutive points are linked, is an oriented graph (according to increasing abscissa), which we will denote by~$ \Gamma_{{\cal W}_m}$.~$V_m$ is called the set of vertices of the graph~$ \Gamma_{{\cal W}_m}$. We will denote, in the sequel, by
$${\cal N}^{\cal S}_m=2\, N_b^m+  N_b-2$$

\noindent the number of vertices of the graph~$ \Gamma_{{\cal W}_m}$, and we will write:

 $$V_m = \left \lbrace {\cal S}_0^m,  {\cal S}_1^m, \hdots,  {\cal S}_{{\cal N}_m-1}^m \right \rbrace $$  %Le graphe~$ \Gamma_{{\cal W}_m}$ est appelé \textbf{sous-graphe} de~$ %\Gamma_{{\cal % W} }$ \\

\end{definition}

 \begin{figure}[h!]
 \center{\psfig{height=10cm,width=12cm,angle=0,file=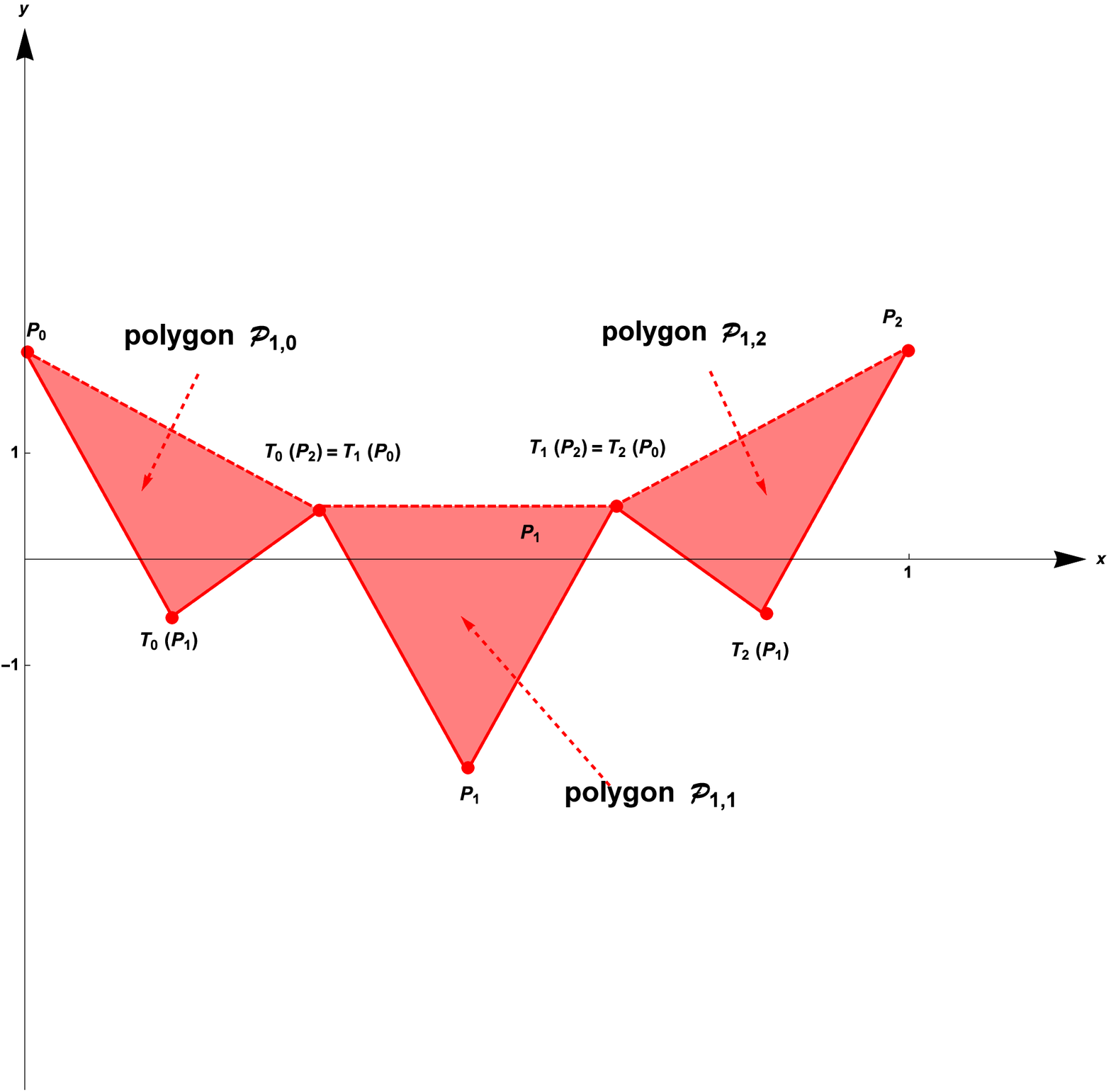}}
\caption{The polygons~${\cal P}_{1,0}$,~${\cal P}_{1,1}$,~${\cal P}_{1,2}$, in the case where~$\lambda= \displaystyle \frac{1}{2}$, and~$N_b=3$.}
 \end{figure}

 \begin{figure}[h!]
 \center{\psfig{height=8cm,width=12cm,angle=0,file=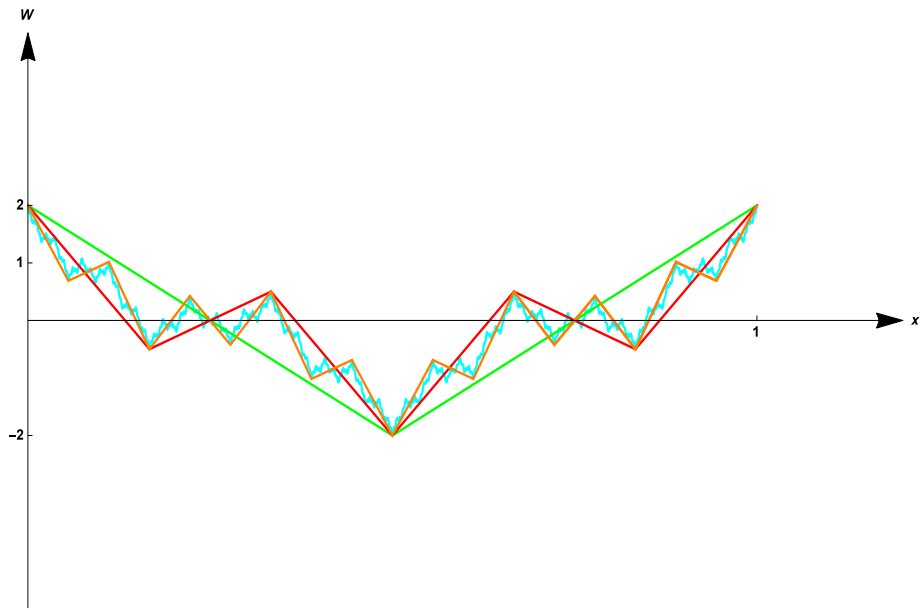}}
\caption{The graphs~$ \Gamma_{{\cal W}_0 }$ (in green), ~$ \Gamma_{{\cal W}_1 }$ (in red),~$ \Gamma_{{\cal W}_2 }$ (in orange),~$ \Gamma_{{\cal W} }$ (in cyan), in the case where~$\lambda= \displaystyle \frac{1}{2}$, and~$N_b=3$.    }
 \end{figure}

\vskip 1cm
\newpage

\begin{definition}\textbf{Consecutive vertices on the graph~$\Gamma_{ \cal W} $ }\\

\noindent Two points~$X$ et~$Y$ de~$\Gamma_{{ \cal W} }$ will be called \textbf{\emph{consecutive vertices}} of the graph~$\Gamma_{ \cal W} $ if there exists a natural integer~$m$, and an integer~$j $ of~\mbox{$\left \lbrace  0,...,N_b-2  \right \rbrace$}, such that:

$$X = \left (T_{i_1}\circ \hdots \circ T_{i_m}\right)(P_j) \quad \text{et} \quad Y = \left (T_{i_1}\circ \hdots \circ T_{i_m}\right)(P_{j+1})
\qquad\left \lbrace i_1,\hdots, i_m \right \rbrace \,\in\,\left \lbrace  0,...,N_b-1  \right \rbrace^m $$

\noindent or:

$$X = \left (T_{i_1}\circ  T_{i_2}\circ \hdots \circ T_{i_m}\right)\left (P_{N_b-1}\right) \quad \text{et} \quad Y =\left (T_{i_1+1}\circ T_{i_2}\hdots \circ T_{i_m} \right)(P_{0})$$

\end{definition}

\vskip 1cm

\begin{definition}
\noindent For any natural integer~$m$, the~$ {\cal N}^{\cal S}_m$ consecutive vertices of the graph~$  \Gamma_{{\cal W}_m} $ are, also, the vertices of~$N_b^m$ simple polygons~${\cal P}_{m,j}$,~\mbox{$0 \leq j \leq N_b^m-1$}, with~$N_b$ sides. For any integer~$j$ such that~\mbox{$0 \leq j \leq N_b^m-1$}, one obtains each polygon by linking the point number~$j$ to the point number~$j+1$ if~\mbox{$j = i \, \text{mod } N_b$},~\mbox{$0 \leq i \leq N_b-2$}, and the point number~$j$ to the point number~$j-N_b+1$ if~\mbox{$j =-1 \, \text{mod } N_b$}. These polygons generate a Borel set of~$\R^2$.

\end{definition}
\vskip 1cm

\begin{definition}\textbf{Polygonal domain delimited by the graph~$  \Gamma_{{\cal W}_m} $,~$m\,\in\,\N $}\\

\noindent For any natural integer~$m$, well call \textbf{polygonal domain delimited by the graph~$  \Gamma_{{\cal W}_m} $}, and denote by~\mbox{$ {\cal D} \left ( \Gamma_{{\cal W}_m}\right) $}, the reunion of the~$N_b^m$ polygons~${\cal P}_{m,j}$,~\mbox{$0 \leq j \leq N_b^m-1$}, with~$N_b$ sides.
\end{definition}
\vskip 1cm

\begin{definition}\textbf{Polygonal domain delimited by the graph~$  \Gamma_{{\cal W} } $ }\\

\noindent We will call \textbf{polygonal domain delimited by the graph~$  \Gamma_{{\cal W}} $}, and denote by~\mbox{$ {\cal D} \left ( \Gamma_{{\cal W} }\right) $}, the limit:
$$ {\cal D} \left ( \Gamma_{{\cal W} }\right)  = \displaystyle \lim_{n \to + \infty} {\cal D} \left ( \Gamma_{{\cal W}_m}\right) $$

\end{definition}

\vskip 1cm
\begin{definition}\textbf{Word, on the graph~$\Gamma_{ \cal W} $}\\

\noindent Let~$m  $ be a strictly positive integer. We will call \textbf{number-letter} any integer~${\cal M}_i$ of~\mbox{$\left \lbrace 0, \hdots, N_b-1 \right \rbrace $}, and \textbf{word of length~$|{\cal M}|=m$}, on the graph~$\Gamma_{ \cal W} $, any set of number-letters of the form:

$${\cal M}=\left ( {\cal M}_1, \hdots, {\cal M}_m\right)$$

\noindent We will write:

$$T_{\cal M}= T_{{\cal M}_1} \circ \hdots \circ  T_{{\cal M}_m}  $$

\end{definition}

\vskip 1cm

\begin{definition}\textbf{Edge relation, on the graph~$\Gamma_{ \cal W} $}\\

\noindent Given a natural integer~$m$, two points~$X$ and~$Y$ of~$\Gamma_{{ \cal W}_m}$ will be called \emph{\textbf{adjacent}} if and only if~$X$ and~$Y$ are two consecutive vertices of~$\Gamma_{{ \cal W}_m}$. We will write:

$$X \underset{m }{\sim}  Y$$

\noindent This edge relation ensures the existence of a word~{${\cal M}=\left ( {\cal M}_1, \hdots, {\cal M}_m\right)$} of length~$ m$, such that~$X$ and~$Y$ both belong to the iterate:

$$T_{\cal M} \,V_0=\left (T_{{\cal M}_1} \circ \hdots \circ  T_{{\cal M}_m} \right) \,V_0$$

 \noindent Given two points~$X$ and~$Y$ of the graph~$\Gamma_{ \cal W} $, we will say that~$X$ and~$Y$ are \textbf{\emph{adjacent}} if and only if there exists a natural integer~$m$ such that:
$$X  \underset{m }{\sim}  Y$$
%\noindent Le graphe~$\Gamma_{ \cal W} $ est dit \emph{\textbf{connecté}}.
\end{definition}

\vskip 1cm

\begin{proposition}\textbf{Adresses, on the graph of the Weierstrass function}\\

\noindent Given a strictly positive integer~$m$, and a word~${\cal M}=\left ( {\cal M}_1, \hdots, {\cal M}_m\right)$ of length~$m\,\in\,\N^\star$, on the graph~$\Gamma_{ {\cal W}_m  }$, for any integer~$j$ of~$\left \lbrace 1,...,N_b-2\right \rbrace  $, any~$X=T_{\cal M}(P_j)$ de~$ V_m \setminus V_{0}$, i.e. distinct from one of the~$N_b $ fixed point~$P_i$,
  ~\mbox{$0 \leq i \leq N_b-1$}, has exactly two adjacent vertices, given by:

$$T_{\cal M}(P_{j+1})\quad \text{et} \quad T_{\cal M}(P_{j-1})$$

\noindent where:

$$T_{\cal M}  = T_{{\cal M}_1} \circ \hdots \circ  T_{{\cal M}_m}   $$

\noindent By convention, the adjacent vertices of~$T_{{\cal M} }(P_{0})  $ are~$T_{{\cal M} }(P_{1})$ and~$T_{{\cal M} }(P_{N_b-1})$,
those of~$T_{{\cal M} }(P_{N_b-1})$,~$T_{{\cal M} }(P_{N_b-2})$ and~$T_{{\cal M} }(P_{0})$ .

\end{proposition}

\vskip 1cm

\begin{definition} \textbf{$m^{th}-$order subcell,~$m \,\in\,\N^\star$, related to a pair of points of the graph~$\Gamma_{\cal W}$}\\

\noindent Given a strictly positive integer~$m$, and two points~$X$ and~$Y$ of~$V_m$ such that~\mbox{$X  \underset{m }{\sim}  Y $}, we will call \textbf{$m^{th}-$order subcell, related to the pair of points~$(X,Y)$}, the polygon, the vertices of which are~$X$,~$Y$, and the intersection points of the edge between the vertices at the extremities of the polygon, i.e. the respective intersection points of polygons of the type~${\cal P}_{m,j-1}$ and~${\cal P}_{m,j }$,~\mbox{$1 \leq j \leq N_b^m-1$}, on the one hand, and of the type~${\cal P}_{m,j }$ and~${\cal P}_{m,j+1 }$,~\mbox{$0 \leq j \leq N_b^m-2$}, on the other hand.
\end{definition}

\vskip 1cm

 \begin{figure}[h!]
 \center{\psfig{height=8cm,width=14cm,angle=0,file=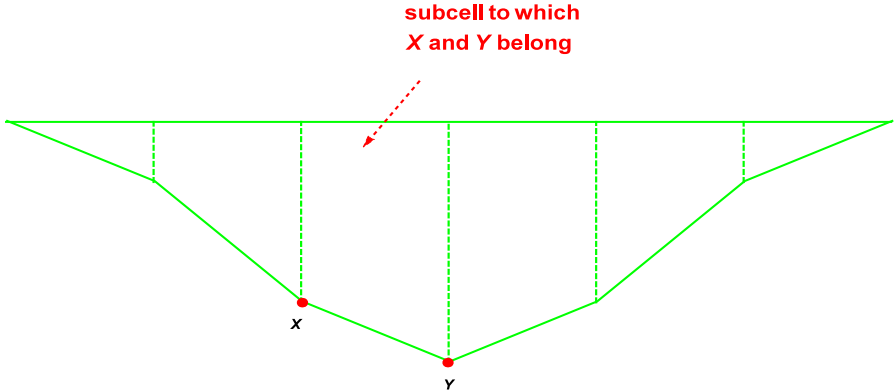}}
\caption{A~$m^{th}-$order subcell, in the case where~$\lambda= \displaystyle \frac{1}{2}$, and~$N_b=7$.}
 \end{figure}

\begin{notation}
\noindent For any integer~$j $ belonging to~\mbox{$\left \lbrace  0,...,N_b-1  \right \rbrace $}, any natural integer~$m$, and any word~$\cal M$ of length~$m$, we set:

 $$\begin{array}{ccc} T_{\cal M}  \left (P_{j } \right) & =&
\left ( x \left (T_{\cal M}  \left (P_{j } \right) \right),  y \left (T_{\cal M}  \left (P_{j } \right) \right) \right)
 \end{array} \quad , \quad \begin{array}{ccc} T_{\cal M}  \left (P_{j+1 } \right) & =&
\left ( x \left (T_{\cal M}  \left (P_{j+1 } \right) \right),  y \left (T_{\cal M}  \left (P_{j+1 } \right) \right) \right)
 \end{array}
 $$

$$L_m= x\left ( T_{\cal M}  \left (P_{j+1 } \right) \right)-x\left ( T_{\cal M}  \left (P_{j } \right)\right)=
  \displaystyle \frac{1}{(N_b-1)\, N_b^m} $$

\end{notation}

\vskip 1cm

\begin{proposition}\textbf{An upper bound and lower bound, for the box-dimension of the graph~$\Gamma_{\cal W}$}\\

\noindent For any integer~$j$ belonging to~\mbox{$\left \lbrace 0, 1, \hdots, N_b-2 \right \rbrace $}, each natural integer~$m$, and each word~$\cal M$ of length~$m$, let us consider the rectangle, the width of which is:

$$L_m= x\left ( T_{\cal M}  \left (P_{j+1 } \right) \right)-x\left ( T_{\cal M}  \left (P_{j } \right)\right)=
  \displaystyle \frac{1}{(N_b-1)\, N_b^m} $$

\noindent and height~$|h_{j,m}|$, such that the points~\mbox{$T_{\cal M}  \left (P_{j+1 }\right) $} and~\mbox{$T_{\cal M}  \left (P_{j+1 } \right )$} are two vertices of this rectangle.\\
     \noindent Then:

  $$L_m^{2-D_{\cal W}}  \,(N_b-1)^{2-D_{\cal W}}\,\left |  \left \lbrace    \displaystyle\frac{2}{1-\lambda}\,   \displaystyle \min_{0 \leq j \leq N_b-1 }  \sin\left ( \displaystyle\frac{ \pi\,(2\,j+1) }{ N_b-1 }\right )    - \displaystyle  \frac{ \pi    }{    N_b\,(N_b-1)\, \left (\lambda\, N_b  -1 \right)}   \right \rbrace  \right| \leq |h_{j,m}|  $$

\noindent and:

  $$  |h_{j,m}| \leq  \eta_{2-D_{\cal W}   }\,L_m^{2-D_{\cal W}}    $$

  \noindent where the real constant~\mbox{$  \eta_{2-D_{\cal W}   }$} is given by :

 $$  \eta_{2-D_{\cal W}   }  = 2\, \pi^2\, (N_b-1)^{2-D_{\cal W}}\, \left \lbrace
   \displaystyle  \frac{   (2\,N_b-1)\, \lambda\, (N_b^2-1) } {(N_b-1)^2 \, (1- \lambda )\,(\lambda \,N_b^{   2}-1) }    +
   \displaystyle  \frac{  2\, N_b } {   (\lambda \,N_b^{ 2  }-1)\, (\lambda \,N_b^{ 3  }-1)  } \right \rbrace
               $$
 
\noindent There exists thus a positive constant
\footnotesize
$$C= \displaystyle \max \,\left \lbrace  \,(N_b-1)^{2-D_{\cal W}}\,\left |  \left \lbrace    \displaystyle\frac{2}{1-\lambda}\,   \displaystyle \min_{0 \leq j \leq N_b-1 }  \sin\left ( \displaystyle\frac{ \pi\,(2\,j+1) }{ N_b-1 }\right )    - \displaystyle  \frac{ \pi    }{    N_b\,(N_b-1)\, \left (\lambda\, N_b  -1 \right)}   \right \rbrace  \right| ,\eta_{2-D_{\cal W}   }  \right \rbrace  $$

\normalsize
\noindent  such that the graph~$\Gamma_{\cal W}$ on~$L_m$ can be covered by at least and at most:

 $$ N_m\, \left \lbrace C\, \left ( \displaystyle \frac{L_m}{N_m} \right)^{1-D_{\cal W}   }+1 \right \rbrace = C\,L_m^{1-D_{\cal W}   }\, N_m^{D_{\cal W}   }+N_m$$

 \noindent squares, the side length of which is~\mbox{$\displaystyle \frac{L_m}{N_m}$}.\\

\end{proposition}

\vskip 1cm

\begin{proof}

For any pair of integers~$(i_m,j)$ of~\mbox{$\left \lbrace  0,...,N_b-2  \right \rbrace^2$}:

 %$$\begin{array}{ccc} T_k\left (P_{j } \right)& =&\left( \displaystyle \frac{x_j+k}{N_b}, \lambda\, y_j + \cos\left(2\,\pi \,\left(\frac{x_j+k}{N_b}\right)\right) \right)
% \end{array}
% $$

 $$\begin{array}{ccc} T_{i_m} \left (P_{j } \right)& =&\left( \displaystyle \frac{x_j+i_m}{N_b}, \lambda\, y_j + \cos\left(2\,\pi \,\left(\frac{x_j+i_m}{N_b}\right)\right) \right)
 \end{array}
 $$

\noindent For any pair of integers~~$(i_m,i_{m-1},j)$ of~\mbox{$\left \lbrace  0,...,N_b-2  \right \rbrace^3$}:

 $$\begin{array}{ccc} T_{i_{m-1}}  \left ( T_{i_m}\left (P_{j } \right)\right)& =&\left( \displaystyle \frac{\frac{x_j+{i_m}}{N_b}+i_{m-1}}{N_b},
  \lambda^2\, y_j + \lambda\,\cos\left(2\,\pi \,\left(  \frac{x_j+{i_m}}{N_b} \right)\right)+ \cos\left(2\,\pi \,\left(\frac{\frac{x_j+{i_m}}{N_b}+i_{m-1} }{N_b}\right)\right) \right)\\ \\
  & =&\left( \displaystyle  \frac{x_j+{i_m}}{ N_b^2}+ \displaystyle \frac{  i_{m-1}}{N_b},
  \lambda^2\, y_j + \lambda\,\cos\left(2\,\pi \,\left(  \frac{x_j+{i_m}}{N_b} \right)\right)+ \cos\left(2\,\pi \,\left( \displaystyle  \frac{x_j+{i_m}}{ N_b^2}+ \displaystyle \frac{  i_{m-1}}{N_b}\right)\right) \right)\\
 \end{array}
 $$
 \noindent For any pair of integers~$(i_m,i_{m-1},i_{m-2},j)$ of~\mbox{$\left \lbrace  0,...,N_b-2  \right \rbrace^4$}:

 $$\begin{array}{ccc} T_{i_{m-2}} \left (T_{i_{m-1}} \left ( T_{i_m}\left (P_{j } \right)\right)\right)& =&
 \bigg(  \displaystyle  \frac{x_j+{i_m}}{ N_b^3}+ \displaystyle \frac{  i_{m-1}}{N_b^2} + \displaystyle  \frac{  i_{m-2}}{N_b},\\
  && \lambda^3\, y_j +  \lambda^2\,  \cos\left(2\,\pi \,\left(  \frac{x_j+ i_m }{N_b} \right)\right) \\
 &&+\lambda \, \cos\left(2\,\pi \,\left( \displaystyle  \frac{x_j+{i_m}}{ N_b^2}+ \displaystyle \frac{  i_{m-1}}{N_b}\right)\right)
  +\cos\left(2\,\pi \,\left( \displaystyle  \frac{x_j+{i_m}}{ N_b^3}+ \displaystyle \frac{  i_{m-1}}{N_b^2} +  \displaystyle  \frac{  i_{m-2}}{N_b}\right)  \right)  \bigg)
 \end{array}
 $$

\noindent Given a strictly positive integer~$m$, and two points~$X$ and~$Y$ of~$V_m$ such that:

$$X  \underset{m }{\sim}  Y $$

\noindent there exists a word~$\cal M$ of length~\mbox{$|{\cal M}|=m$}, on the graph~$\Gamma_{ \cal W} $, and an integer~$j$ of~\mbox{$\left \lbrace  0,...,N_b-2  \right \rbrace^2$}, such that:

$$ X= T_{\cal M}  \left (P_{j } \right) \quad   , \quad Y= T_{\cal M}  \left (P_{j+1 } \right)$$

\noindent Let us write~$T_{\cal M}$ under the form:

$$T_{\cal M}=T_{i_m} \circ T_{i_{m-1}} \circ \hdots \circ T_{i_1}$$

\noindent where~\mbox{$(i_1,\hdots,i_m)\,\in\, \left \lbrace  0,...,N_b-1  \right \rbrace^m$}.\\

\noindent One has then:

 $$  x \left (T_{\cal M}  \left (P_{j } \right) \right )  =
    \displaystyle  \frac{x_j }{ N_b^m}+ \displaystyle \sum_{k=1}^m\frac{  i_{k}}{N_b^k}  \quad , \quad
       x \left (T_{\cal M}  \left (P_{j+1 } \right) \right) =
    \displaystyle  \frac{x_{j+1} }{ N_b^m}+ \displaystyle \sum_{k=1}^m\frac{  i_{k}}{N_b^k}
 $$

\noindent and:
  \footnotesize

 $$\left \lbrace \begin{array}{ccc}  y\left (T_{\cal M}  \left (P_{j } \right) \right) & =&\lambda^m\, y_{j } + \displaystyle \sum_{k=1}^{m } \lambda^{m-k}\,  \cos\left(2\,\pi \,\left(   \displaystyle  \frac{x_{j  } }{ N_b^{k}}+
    \displaystyle \sum_{\ell=0}^{k}\frac{  i_{ m-\ell}}{N_b^{k- \ell}}\right)\right)
        \\ \\
       y \left (T_{\cal M}  \left (P_{j+1 } \right) \right) & =&\lambda^m\, y_{j+1 } + \displaystyle \sum_{k=1}^{m } \lambda^{m-k}\,  \cos\left(2\,\pi \,\left(   \displaystyle  \frac{x_{j+1 } }{ N_b^{k}}+
    \displaystyle \sum_{\ell=0}^{k}\frac{  i_{ m-\ell}}{N_b^{k- \ell}}\right)\right)
 \end{array} \right.
 $$

\normalsize

\noindent This leads to:

 $$\begin{array}{ccc}    h_{j,m} -\lambda^m\, \left (y_{j+1 }-y_{j } \right) & =  &
            \displaystyle \sum_{k=1}^{m } \lambda^{m-k}\, \left \lbrace  \cos\left(2\,\pi \,\left(   \displaystyle  \frac{x_{j+1 } }{ N_b^{k}}+
    \displaystyle \sum_{\ell=0}^{k}\frac{  i_{ m-\ell}}{N_b^{k- \ell}}\right)\right) -\cos\left(2\,\pi \,\left(   \displaystyle  \frac{x_{j  } }{ N_b^{k}}  -
  \displaystyle \sum_{\ell=0}^{k}\frac{  i_{ m-\ell}}{N_b^{k- \ell}}\right)\right) \right \rbrace \\
    \\ \\
  & =  & -2 \displaystyle \sum_{k=1}^{m } \lambda^{m-k}\, \sin \left( \pi \,\left(   \displaystyle  \frac{x_{j+1 }-x_j }{  N_b^{k}} \right) \right) \, \sin \left(2\,\pi \,\left(   \displaystyle  \frac{x_{j+1 }+x_j }{ 2\,N_b^{k}}+
    \displaystyle \sum_{\ell=0}^{k}\frac{  i_{ m-\ell}}{N_b^{k- \ell}}\right) \right)  \\
    \\ \\
         \end{array}
 $$

\noindent Taking into account:

 $$\begin{array}{ccc}   \lambda^m\, \left (y_{j+1 }-y_{j } \right) & =  &
         \displaystyle\frac{\lambda^m}{1-\lambda}\, \left ( \cos\left ( \displaystyle\frac{2\,\pi\,(j+1)}{N_b-1}\right )- \cos\left ( \displaystyle\frac{2\,\pi\,j}{N_b-1}\right ) \right) \\
%  & =  &         -2\,  \displaystyle\frac{\lambda^m}{1-\lambda}\,   \sin\left ( \displaystyle\frac{2\,\pi\,(j+1-j)}{2\,(N_b-1)}\right )\, \sin\left ( %\displaystyle\frac{2\,\pi\,(j+1+j)}{2\, (N_b-1)}\right ) \\
   & =  &         -2\,  \displaystyle\frac{\lambda^m}{1-\lambda}\,   \sin\left ( \displaystyle\frac{ \pi }{ N_b-1 }\right )\, \sin\left ( \displaystyle\frac{ \pi\,(2\,j+1 )}{  N_b-1 }\right ) \\

         \end{array}
 $$

\noindent one has:

\footnotesize

 $$\begin{array}{ccc}    h_{j,m}+ 2\,  \displaystyle\frac{\lambda^m}{1-\lambda}\,   \sin\left ( \displaystyle\frac{ \pi }{ N_b-1 }\right )\, \sin\left ( \displaystyle\frac{ \pi\,(2\,j+1 )}{  N_b-1 }\right )
  & =  &  -2\,\displaystyle \sum_{k=1}^{m } \lambda^{m-k}\, \sin \left(   \displaystyle  \frac{ \pi   }{  N_b^{k+1} \,(N_b-1)} \right)   \, \sin \left(   \displaystyle  \frac{ \pi \,(2\,j+1) }{  N_b^{k+1}\,(N_b-1)}+2\, \pi \,
    \displaystyle \sum_{\ell=0}^{k}\frac{  i_{ m-\ell}}{N_b^{k- \ell}}\right)  \\
         \end{array}
 $$

\normalsize

%\noindent i.e.:

 %$$\begin{array}{ccc}    y \left (T_{\cal M}  \left (P_{j  } \right) \right)- y \left (T_{\cal M}  \left (P_{j+1  } \right) \right)
%  & =  &
 %        2\,   \displaystyle\frac{\lambda^m}{1-\lambda}\,   \sin\left ( \displaystyle\frac{ \pi }{ N_b-1 }\right )\, \sin\left ( \displaystyle\frac{ \pi\,(2\,j+1 )}{  N_b-1 % % % % }\right ) \\
 %        &&+ 2\,\displaystyle \sum_{k=1}^{m } \lambda^{m-k}\, \sin \left(   \displaystyle  \frac{ \pi   }{  N_b^{k+1} \,(N_b-1)} \right)   \, \sin \left(   \displaystyle  \frac{ % \pi \,(2\,j+1) }{  N_b^{k+1}\,(N_b-1)}+2\, \pi \,
   %   \displaystyle \sum_{\ell=0}^{k}\frac{  i_{ m-\ell}}{N_b^{k- \ell}}\right)  \\
   %   \\ \\
   %       \end{array}
  %$$

\noindent Thus:

 $$\begin{array}{ccc}  \left |  y \left (T_{\cal M}  \left (P_{j  } \right) \right)- y \left (T_{\cal M}  \left (P_{j+1  } \right) \right)
 - 2\,   \displaystyle\frac{\lambda^m}{1-\lambda}\,   \sin\left ( \displaystyle\frac{ \pi }{ N_b-1 }\right )\, \sin\left ( \displaystyle\frac{ \pi\,(2\,j+1 )}{  N_b-1 }\right ) \right |
     & \leq  &   \displaystyle \sum_{k=1}^{m }     \displaystyle  \frac{ 2\,\pi \,\lambda^{m-k}  }{  N_b^{k+1} \,(N_b-1)}      \\
   & = &    \displaystyle  \frac{ \pi \,\lambda^{m  } \, \left (1-\frac{1}{\lambda^m\, N_b^m } \right) }{ \lambda\, N_b  \,N_b\,(N_b-1)\, \left (1-\frac{1}{\lambda \, N_b  } \right)}      \\
    & \leq &    \displaystyle  \frac{ \pi   \, \lambda^{m } }{    N_b\,(N_b-1)\, \left (\lambda\, N_b  -1 \right)}      \\
      \\ \\
         \end{array}
 $$

\normalsize

\noindent which leads to:

 $$\begin{array}{ccc}     y \left (T_{\cal M}  \left (P_{j  } \right) \right)- y \left (T_{\cal M}  \left (P_{j+1  } \right) \right)
     & \geq  & 2\,   \displaystyle\frac{\lambda^m}{1-\lambda}\,   \sin\left ( \displaystyle\frac{ \pi }{ N_b-1 }\right )\, \sin\left ( \displaystyle\frac{ \pi\,(2\,j+1 )}{  N_b-1 }\right )    - \displaystyle  \frac{ \pi   \, \lambda^{m } }{    N_b\,(N_b-1)\, \left (\lambda\, N_b  -1 \right)}      \\
      \\ \\
         \end{array}
 $$

\noindent or:

 $$\begin{array}{ccc}     y \left (T_{\cal M}  \left (P_{j+1  } \right) \right)- y \left (T_{\cal M}  \left (P_{j   } \right) \right)
     & \geq  & 2\,   \displaystyle\frac{\lambda^m}{1-\lambda}\,   \sin\left ( \displaystyle\frac{ \pi }{ N_b-1 }\right )\, \sin\left ( \displaystyle\frac{ \pi\,(2\,j+1 )}{  N_b-1 }\right )    - \displaystyle  \frac{ \pi   \, \lambda^{m } }{    N_b\,(N_b-1)\, \left (\lambda\, N_b  -1 \right)}      \\
      \\ \\
         \end{array}
 $$

\noindent Due to the symmetric roles played by~$T_{\cal M}  \left (P_{j   } \right )$ and~$T_{\cal M}  \left (P_{j+1  } \right )$, one may only consider the case when:

 $$\begin{array}{ccccc}     y \left (T_{\cal M}  \left (P_{j  } \right) \right)- y \left (T_{\cal M}  \left (P_{j+1  } \right) \right)
     & \geq  & 2\,   \displaystyle\frac{\lambda^m}{1-\lambda}\,   \sin\left ( \displaystyle\frac{ \pi }{ N_b-1 }\right )\, \sin\left ( \displaystyle\frac{ \pi\,(2\,j+1 )}{  N_b-1 }\right )    - \displaystyle  \frac{ \pi   \, \lambda^{m } }{    N_b\,(N_b-1)\, \left (\lambda\, N_b  -1 \right)}    & \geq & 0   \\
      & \geq  &  \lambda^m\, \left \lbrace    \displaystyle\frac{2}{1-\lambda}\,   \displaystyle \min_{0 \leq j \leq N_b-1 }  \sin\left ( \displaystyle\frac{ \pi\,(2\,j+1) }{ N_b-1 }\right )    - \displaystyle  \frac{ \pi    }{    N_b\,(N_b-1)\, \left (\lambda\, N_b  -1 \right)}   \right \rbrace     \\
      \\ \\
         \end{array}
 $$

%\noindent Pour~$\ell=k$ :

%$$ 2\,\pi\,    i_{ m-k}   \,\in\,2\,\pi\,\N$$

\noindent The predominant term is thus:

$$   \lambda^{m  }
=  e^{ m\,(D_{\cal W}-2)\, \ln N_b}= N_b^{ m\,(D_{\cal W}-2)}  = L_m^{2-D_{\cal W}}  \,(N_b-1)^{2-D_{\cal W}}   $$
%\noindent or:

%$$ \displaystyle  \frac{1}{N_b^{2m}}  \quad \text{if}  \quad  \displaystyle  \frac{1}{N_b^{2 }} > \lambda  $$

 %$$D_{\cal W}= 2+\displaystyle \frac{\ln \lambda}{\ln N_b} \quad , \quad \lambda= e^{(D_{\cal W}-2)\, \ln N_b}= N_b^{(D_{\cal W}-4) } $$

% \noindent Ainsi :

% $$ N_b^{ m\,(D_{\cal W}-2)} \lesssim h_{j,m} \lesssim  N_b^{m\,(D_{\cal W}-2)}   $$

\noindent One also has:

 $$\begin{array}{ccc}  \left | h_{j,m}\right| & \leq  &
          \displaystyle\frac{2\,\lambda^m }{1-\lambda}\,  \displaystyle\frac{ \pi^2\,(2\,j+1)}{(N_b-1)^2}   +
  2\,\displaystyle \sum_{k=1}^{m } \lambda^{m-k}\, \pi \,\left\lbrace   \displaystyle \frac{2\,j+1}{(N_b-1)\, N_b^k}+
  2\, \displaystyle \sum_{\ell=0}^{k}\frac{  i_{ m-\ell}}{N_b^{k- \ell}}\right \rbrace \,   \displaystyle  \frac{  \pi }{ (N_b-1)\, N_b^k}
              \\ \\
         & =&       \displaystyle\frac{2\,\lambda^m}{1-\lambda}\,  \displaystyle\frac{ \pi^2\,(2\,j+1)}{(N_b-1)^2}  +
   \displaystyle  \frac{  2\,\pi^2\, \lambda^m  }{  N_b-1 }\, \displaystyle \sum_{k=1}^{m }  \,\left\lbrace   \displaystyle \frac{(2\,j+1)\, \lambda^{ -k}}{(N_b-1)\, N_b^{2k}}+
  2\, \displaystyle \sum_{\ell=0}^{k}\frac{  i_{ m-\ell}\,  \lambda^{ -k} }{N_b^{2k- \ell}}\right \rbrace
         \\ \\
        & = &      \displaystyle\frac{2\,\lambda^m}{1-\lambda}\,  \displaystyle\frac{ \pi^2\,(2\,j+1)}{(N_b-1)^2} \\
 && +
   \displaystyle  \frac{  2\,\pi^2\, \lambda^m }{  N_b-1 }\,   \left\lbrace
    \displaystyle \frac{\lambda^{-1}\,N_b^{-2}\,(2\,j+1) }{(N_b-1) } \,  \displaystyle \frac{  (1- \lambda^{ -m }\,N_b^{ -2m })}{1- \lambda^{ - 1}\,N_b^{ - 2}} +
  2\, \displaystyle \sum_{k=1}^m \frac{  (N_b-1)\,  \lambda^{ -k}
  }{N_b^{2k }} \displaystyle \frac{1-N_b^{-k-1}}{1-N_b^{-1}}\right \rbrace \\ \\
    & \leq  &      \displaystyle\frac{2\,\lambda^m}{1-\lambda}\,  \displaystyle\frac{ \pi^2\,(2\,N_b-1)}{(N_b-1)^2}
   +
   \displaystyle  \frac{  2\,\pi^2\, \lambda^m }{  N_b-1 }\,
    \displaystyle \frac{ (2\,N_b-1) }{(N_b-1) } \,  \displaystyle \frac{  (1- \lambda^{ -m }\,N_b^{ -2m })}{  \lambda \,N_b^{   2}-1}   \\ \\
 && +
   \displaystyle  \frac{  2\,\pi^2\, \lambda^m }{  N_b-1 }\,
  2\, \displaystyle  \frac{ \lambda^{ -1 }\,N_b^{ -2  }\, (N_b-1)\, (1- \lambda^{ -m }\,N_b^{ -2m })
  }{ (1-N_b^{-1})\,(1-\lambda^{ -1 }\,N_b^{ -2  })}  \\ \\
 && -
   \displaystyle  \frac{  2\,\pi^2\, \lambda^m }{  N_b-1 }\,
  2\, \displaystyle  \frac{ \lambda^{ -1 }\,N_b^{ -3  }\, (N_b-1)\, (1- \lambda^{ -m }\,N_b^{ -3m })}{ (1-N_b^{-1})\,(1-\lambda^{ -1 }\,N_b^{ -3  })}
  \\ \\
 %    & =  &      \displaystyle\frac{2\,\lambda^m}{1-\lambda}\,  \displaystyle\frac{ \pi^2\,(2\,N_b-1)}{(N_b-1)^2} \\
 %&& +
  % \displaystyle  \frac{  2\,\pi^2\, \lambda^m }{  N_b-1 }\,
  %  \displaystyle \frac{ (2\,N_b-1) }{(N_b-1) } \,  \displaystyle \frac{  (1- \lambda^{ -m }\,N_b^{ -2m })}{  \lambda \,N_b^{   2}-1}   \\ \\
 %&& +
 %  \displaystyle  \frac{  4\,\pi^2\, \lambda^m }{  N_b-1 }\,
 %   \displaystyle  \frac{  N_b\, (1- \lambda^{ -m }\,N_b^{ -2m })
 % }{   \lambda^{   }\,N_b^{  2  }-1 }  \\ \\
%&& -
 %  \displaystyle  \frac{  4\,\pi^2\, \lambda^m }{  N_b-1 }\,
  %  \displaystyle  \frac{   N_b \, (1- \lambda^{ -m }\,N_b^{ -3m })}{   \lambda \,N_b^{  3  }-1 }
 % \\ \\
   & \leq  &      \displaystyle\frac{2\,\lambda^m}{1-\lambda}\,  \displaystyle\frac{ \pi^2\,(2\,N_b-1)}{(N_b-1)^2}  +
   \displaystyle  \frac{  2\,\pi^2\, \lambda^m }{  N_b-1 }\,
    \displaystyle \frac{ (2\,N_b-1) }{(N_b-1) } \,  \displaystyle \frac{ 1}{  \lambda \,N_b^{   2}-1}   \\ \\
 && +
   \displaystyle  \frac{  4\,\pi^2\, N_b\, \lambda^m }{   N_b-1  }\,\left \lbrace   \displaystyle  \frac{ 1}{   \lambda \,N_b^{  2  }-1  }  -
   \displaystyle  \frac{ 1}{   \lambda \,N_b^{  3  }-1  }\right \rbrace
  \\ \\
  & =  &2\, \pi^2\,\lambda^m \,\left \lbrace
   \displaystyle  \frac{   (2\,N_b-1)\, \lambda\, (N_b^2-1) } {(N_b-1)^2 \, (1- \lambda )\,(\lambda \,N_b^{   2}-1) }    +
   \displaystyle  \frac{  2\, N_b } {   (\lambda \,N_b^{ 2  }-1)\, (\lambda \,N_b^{ 3  }-1)  } \right \rbrace
  \\ \\

 \end{array}
 $$

 \noindent Since:

$$ x\left ( T_{\cal M}  \left (P_{j+1 } \right) \right)-x\left ( T_{\cal M}  \left (P_{j } \right)\right)=
  \displaystyle \frac{1}{(N_b-1)\, N_b^m} $$

 \noindent and:

 $$D_{\cal W}= 2+\displaystyle \frac{\ln \lambda}{\ln N_b} \quad , \quad \lambda= e^{(D_{\cal W}-2)\, \ln N_b}= N_b^{(D_{\cal W}-2) } $$

 \noindent one has thus:

 $$\begin{array}{ccc}  \left | h_{j,m}\right| & \leq  & 2\, \pi^2\,L_m^{2-D_{\cal W}   }\, \left (N_b-1\right)^{2-D_{\cal W}   }\,\,\left \lbrace
   \displaystyle  \frac{   (2\,N_b-1)\, \lambda\, (N_b^2-1) } {(N_b-1)^2 \, (1- \lambda )\,(\lambda \,N_b^{   2}-1) }    +
   \displaystyle  \frac{  2\, N_b } {   (\lambda \,N_b^{ 2  }-1)\, (\lambda \,N_b^{ 3  }-1)  } \right \rbrace
  \\ \\
 \end{array}
 $$

 \vskip 1cm

\end{proof}
\vskip 1cm

\begin{pte}\textbf{Exact computation of the measure of the surface of a simple polygon~${\cal P}_{m,j}$,\\~\mbox{$0 \leq j \leq N_b^m-1$}, with~$N_b$ sides} \\

\noindent Let us note that, given a natural integer~$m$, there exists a word~$\cal M$ of length~$m$ such that the ordered set, according to increasing abscissa, of the~$N_b$ vertices of a simple polygon~${\cal P}_{m,j}$,~\mbox{$0 \leq j \leq N_b^m-1$}, can be written as:

$$T_{\cal M}\left (P_0\right) \quad , \quad T_{\cal M}\left (P_1\right)\quad , \quad \hdots \quad , \quad T_{\cal M}\left (P_{N_b-1}\right)$$

\noindent This enables one to exactly compute the measure, with respect to the standard Lebesgue measure on~$\R^2$, of any of the aforementioned polygons, as:

\begin{enumerate}
\item[ {i}.] In the case where~$N_b=3$:

$$\begin{array}{ccc}{\cal A}  \left ( {\cal P}_{m,j} \right)&=&
\displaystyle \frac{1}{2}  \left \| \overrightarrow{T_{\cal M}\left (P_0\right)\, T_{\cal M}\left (P_1\right)}\wedge \overrightarrow{T_{\cal M}\left (P_0\right)\, T_{\cal M}\left (P_{2}\right)} \right\|    \\
\end{array} $$

\item[ {ii}.] In the case where~$N_b\geq 4$:

$$\begin{array}{ccc}{\cal A}  \left ( {\cal P}_{m,j} \right)&=&
\displaystyle \frac{1}{2} \, \displaystyle \sum_{j=0}^{N_b-2} \left \| \overrightarrow{T_{\cal M}\left (P_j\right)\, T_{\cal M}\left (P_{j+1}\right)}\wedge \overrightarrow{T_{\cal M}\left (P_j\right)\, T_{\cal M}\left (P_{N_b-1}\right)} \right\| \\
\end{array} $$

\end{enumerate}

\end{pte}

\vskip 1cm

\begin{remark}

\noindent One obtains, for~\mbox{$0 \leq j \leq N_b^m-1$}:

$$\begin{array}{ccc}
&&\left \| \overrightarrow{T_{\cal M}\left (P_j\right)\, T_{\cal M}\left (P_{j+1}\right)}\wedge \overrightarrow{T_{\cal M}\left (P_k\right)\, T_{\cal M}\left (P_{N_b-1}\right)} \right\|=\\ \\&&=\left |
x \left ( \overrightarrow{T_{\cal M}\left (P_j\right)\, T_{\cal M}\left (P_{j+1}\right)}\right)\, y\left ( \overrightarrow{T_{\cal M}\left (P_j\right)\, T_{\cal M}\left (P_{N_b-1}\right)} \right)
-y \left ( \overrightarrow{T_{\cal M}\left (P_j\right)\, T_{\cal M}\left (P_{j+1}\right)}\right)\, x\left ( \overrightarrow{T_{\cal M}\left (P_j\right)\, T_{\cal M}\left (P_{N_b-1}\right)} \right)\right|
\end{array} $$

%$$\begin{array}{ccc} ¨%x \left ( \overrightarrow{T_{\cal M}\left (P_j\right)\, T_{\cal M}\left (P_{j+1}\right)}\right)&=&
%x \left (   T_{\cal M}\left (P_{j+1} \right)\right)-x \left (   T_{\cal M}\left (P_{j } \right)\right)\\
%%&=&    \displaystyle  \frac{1 }{(N_b-1)\, N_b^m}
%\end{array} $$

\noindent where:

$$\begin{array}{ccc}
x \left ( \overrightarrow{T_{\cal M}\left (P_j\right)\, T_{\cal M}\left (P_{N_b-1}\right)}\right)&=&
x \left (   T_{\cal M}\left (P_{j+1} \right)\right)-x \left (   T_{\cal M}\left (P_{N_b-1 } \right)\right)\\
&=&    \displaystyle  \frac{x_j }{ N_b^m}+ \displaystyle \sum_{k=1}^m\frac{  i_{k}}{N_b^k}
- \displaystyle  \frac{x_{N_b-1} }{ N_b^m}- \displaystyle \sum_{k=1}^m\frac{  i_{k}}{N_b^k} \\
&=&    \displaystyle  \frac{N_b-1-j }{(N_b-1)\, N_b^m}
\end{array} $$

\noindent and:

$$\begin{array}{ccc}
y \left ( \overrightarrow{T_{\cal M}\left (P_j\right)\, T_{\cal M}\left (P_{j+1}\right)}\right)&=&
y \left (   T_{\cal M}\left (P_{j+1} \right)\right)-y \left (   T_{\cal M}\left (P_{j } \right)\right)\\ \\
&=&    \lambda^m\, y_{j+1 } + \displaystyle \sum_{k=1}^{m-1} \lambda^{m-k}\,  \cos\left(2\,\pi \,\left(   \displaystyle  \frac{x_{j+1 } }{ N_b^{k}}+
    \displaystyle \sum_{\ell=0}^{k}\frac{  i_{ m-\ell}}{N_b^{k- \ell}}\right)\right)
  +\cos\left(2\,\pi \,\left(  \displaystyle  \frac{x_{j+1 } }{ N_b^m}+ \displaystyle \sum_{k=1}^m\frac{  i_{k}}{N_b^k} \right)  \right) \\
  &&-\lambda^m\, y_{j } + \displaystyle \sum_{k=1}^{m-1} \lambda^{m-k}\,  \cos\left(2\,\pi \,\left(   \displaystyle  \frac{x_{j  } }{ N_b^{k}}+
    \displaystyle \sum_{\ell=0}^{k}\frac{  i_{ m-\ell}}{N_b^{k- \ell}}\right)\right)
  +\cos\left(2\,\pi \,\left(  \displaystyle  \frac{x_{j  } }{ N_b^m}+ \displaystyle \sum_{k=1}^m\frac{  i_{k}}{N_b^k} \right)  \right)\\ \\
  &=&    \lambda^m\,\left ( y_{j+1 }-y_j\right) \\
  &&+ \displaystyle \sum_{k=1}^{m-1} \lambda^{m-k}\,  \left \lbrace \cos\left(2\,\pi \,\left(   \displaystyle  \frac{x_{j+1 } }{ N_b^{k}}+
    \displaystyle \sum_{\ell=0}^{k}\frac{  i_{ m-\ell}}{N_b^{k- \ell}}\right)\right)
    -\cos\left(2\,\pi \,\left(   \displaystyle  \frac{x_{j  } }{ N_b^{k}}+
    \displaystyle \sum_{\ell=0}^{k}\frac{  i_{ m-\ell}}{N_b^{k- \ell}}\right)\right) \right \rbrace \\
    &&+ \cos\left(2\,\pi \,\left(  \displaystyle  \frac{x_{j+1 } }{ N_b^m}+ \displaystyle \sum_{k=1}^m\frac{  i_{k}}{N_b^k} \right)  \right)
  - \cos\left(2\,\pi \,\left(  \displaystyle  \frac{x_{j  } }{ N_b^m}+ \displaystyle \sum_{k=1}^m\frac{  i_{k}}{N_b^k} \right)  \right) \\
 \\
 &=&    \displaystyle\frac{\lambda^m}{1-\lambda}\,\left \lbrace \cos\left ( \displaystyle\frac{2\,\pi\,(j+1)}{N_b-1}\right )
-\cos\left ( \displaystyle\frac{2\,\pi\,j}{N_b-1}\right ) \right \rbrace \\
  &&+ \displaystyle \sum_{k=1}^{m-1} \lambda^{m-k}\,  \left \lbrace \cos\left(2\,\pi \,\left(   \displaystyle  \frac{x_{j+1 } }{ N_b^{k}}+
    \displaystyle \sum_{\ell=0}^{k}\frac{  i_{ m-\ell}}{N_b^{k- \ell}}\right)\right)
    -\cos\left(2\,\pi \,\left(   \displaystyle  \frac{x_{j  } }{ N_b^{k}}+
    \displaystyle \sum_{\ell=0}^{k}\frac{  i_{ m-\ell}}{N_b^{k- \ell}}\right)\right) \right \rbrace \\
    &&+ \cos\left(2\,\pi \,\left(  \displaystyle  \frac{x_{j+1 } }{ N_b^m}+ \displaystyle \sum_{k=1}^m\frac{  i_{k}}{N_b^k} \right)  \right)
  - \cos\left(2\,\pi \,\left(  \displaystyle  \frac{x_{j  } }{ N_b^m}+ \displaystyle \sum_{k=1}^m\frac{  i_{k}}{N_b^k} \right)  \right) \\
 \\
\end{array} $$

 %$$xy_j=\displaystyle\frac{j}{ N_b-1} $$

 %$$y_j=\displaystyle\frac{1}{1-\lambda}\,\cos\left ( \displaystyle\frac{2\,\pi\,j}{N_b-1}\right )$$

 %$$y_{j+1}-y_j=\displaystyle\frac{1}{1-\lambda}\,\left \lbrace \cos\left ( \displaystyle\frac{2\,\pi\,(j+1)}{N_b-1}\right )
 %-\cos\left ( \displaystyle\frac{2\,\pi\,j}{N_b-1}\right ) \right \rbrace $$

 %$$y_{N_b-1}-y_j=\displaystyle\frac{1}{1-\lambda}\,\left \lbrace \cos\left ( \displaystyle\frac{2\,\pi\,(N_b-1)}{N_b-1}\right )
 %-\cos\left ( \displaystyle\frac{2\,\pi\,j}{N_b-1}\right ) \right \rbrace
 %=\displaystyle\frac{1}{1-\lambda}\,\left \lbrace 1
 %-\cos\left ( \displaystyle\frac{2\,\pi\,j}{N_b-1}\right ) \right \rbrace
 %$$

$$\begin{array}{ccc}
y \left ( \overrightarrow{T_{\cal M}\left (P_j\right)\, T_{\cal M}\left (P_{N_b-1}\right)}\right)&=&
y \left (   T_{\cal M}\left (P_{N_b-1} \right)\right)-y \left (   T_{\cal M}\left (P_{j } \right)\right)\\ \\
   &=&    \lambda^m\,\left ( y_{N_b-1}-y_j\right) \\
  &&+ \displaystyle \sum_{k=1}^{m-1} \lambda^{m-k}\,  \left \lbrace \cos\left(2\,\pi \,\left(   \displaystyle  \frac{x_{N_b-1} }{ N_b^{k}}+
    \displaystyle \sum_{\ell=0}^{k}\frac{  i_{ m-\ell}}{N_b^{k- \ell}}\right)\right)
    -\cos\left(2\,\pi \,\left(   \displaystyle  \frac{x_{j  } }{ N_b^{k}}+
    \displaystyle \sum_{\ell=0}^{k}\frac{  i_{ m-\ell}}{N_b^{k- \ell}}\right)\right) \right \rbrace \\
    &&+ \cos\left(2\,\pi \,\left(  \displaystyle  \frac{x_{N_b-1} }{ N_b^m}+ \displaystyle \sum_{k=1}^m\frac{  i_{k}}{N_b^k} \right)  \right)
  - \cos\left(2\,\pi \,\left(  \displaystyle  \frac{x_{j  } }{ N_b^m}+ \displaystyle \sum_{k=1}^m\frac{  i_{k}}{N_b^k} \right)  \right) \\
 \\
\end{array} $$

\noindent Thus:

$$\begin{array}{ccc}
&(N_b-1)\, N_b^m \,\left \| \overrightarrow{T_{\cal M}\left (P_j\right)\, T_{\cal M}\left (P_{j+1}\right)}\wedge \overrightarrow{T_{\cal M}\left (P_k\right)\, T_{\cal M}\left (P_{N_b-1}\right)} \right\| = \\ \\
& \biggl |   \displaystyle\frac{\lambda^m}{1-\lambda}\,\left \lbrace 1
-\cos\left ( \displaystyle\frac{2\,\pi\,j}{N_b-1}\right ) \right \rbrace \\
  & + \displaystyle \sum_{k=1}^{m-1} \lambda^{m-k}\,  \left \lbrace \cos\left(2\,\pi \,\left(   \displaystyle  \frac{1 }{ N_b^{k}}+
    \displaystyle \sum_{\ell=0}^{k}\frac{  i_{ m-\ell}}{N_b^{k- \ell}}\right)\right)
    -\cos\left(2\,\pi \,\left(   \displaystyle  \frac{j }{N_b-1)\, N_b^{k}}+
    \displaystyle \sum_{\ell=0}^{k}\frac{  i_{ m-\ell}}{N_b^{k- \ell}}\right)\right) \right \rbrace \\
    & + \cos\left(2\,\pi \,\left(  \displaystyle  \frac{1 }{ N_b^m}+ \displaystyle \sum_{k=1}^m\frac{  i_{k}}{N_b^k} \right)  \right)
  - \cos\left(2\,\pi \,\left(  \displaystyle  \frac{j}{ N_b-1)\,N_b^m}+ \displaystyle \sum_{k=1}^m\frac{  i_{k}}{N_b^k} \right)  \right) \\
   &+  \displaystyle\frac{(N_b-1-j)\, \lambda^m}{1-\lambda}\,\left \lbrace \cos\left ( \displaystyle\frac{2\,\pi\,(j+1)}{N_b-1}\right )
-\cos\left ( \displaystyle\frac{2\,\pi\,j}{N_b-1}\right ) \right \rbrace  \\
  & + (N_b-1-j)\,\left \lbrace \displaystyle \sum_{k=1}^{m-1} \lambda^{m-k}\,  \left \lbrace \cos\left(2\,\pi \,\left(   \displaystyle  \frac{(j+1) }{N_b-1)\, N_b^{k}}+
    \displaystyle \sum_{\ell=0}^{k}\frac{  i_{ m-\ell}}{N_b^{k- \ell}}\right)\right)
    -\cos\left(2\,\pi \,\left(   \displaystyle  \frac{j }{N_b-1)\, N_b^{k}}+
    \displaystyle \sum_{\ell=0}^{k}\frac{  i_{ m-\ell}}{N_b^{k- \ell}}\right)\right) \right \rbrace  \right \rbrace \\
    & +(N_b-1-j)\,\left \lbrace  \cos\left(2\,\pi \,\left(  \displaystyle  \frac{(j+1) }{N_b-1)\, N_b^m}+ \displaystyle \sum_{k=1}^m\frac{  i_{k}}{N_b^k} \right)  \right)
  - \cos\left(2\,\pi \,\left(  \displaystyle  \frac{j }{(N_b-1)\, N_b^m}+ \displaystyle \sum_{k=1}^m\frac{  i_{k}}{N_b^k} \right)  \right) \right \rbrace
\biggl|
\end{array} $$

\end{remark}

\normalsize

\vskip 1cm

\begin{definition}\textbf{Measure, on the domain delimited by the graph~$  \Gamma_{{\cal W} } $ }\\

\noindent We will call \textbf{domain delimited by the graph~~$  \Gamma_{{\cal W}} $}, and denote by~\mbox{$ {\cal D} \left ( \Gamma_{{\cal W} }\right) $}, the limit:
$$ {\cal D} \left ( \Gamma_{{\cal W} }\right)  = \displaystyle \lim_{n \to + \infty} {\cal D} \left ( \Gamma_{{\cal W}_m}\right) $$

\noindent which has to be understood in the following way: given a continuous function~$u$ on the graph~$\Gamma_{\cal W}$, and a measure with full support~$\mu$ on~$\R^2$, then:

$$\displaystyle \int_{ {\cal D} \left ( \Gamma_{{\cal W} }\right)} u\,d\mu  = \displaystyle \lim_{m \to + \infty}
\displaystyle \sum_{j=0}^{N_b^m-1}  \displaystyle \sum_{X \, \text{vertex of }{\cal P}_{m,j} }u\left ( X \right) \,\mu \left (  {\cal P}_{m,j}  \right)$$

\noindent We will say that~$\mu$ is a \textbf{measure, on the domain delimited by the graph~$  \Gamma_{{\cal W} } $}.
\end{definition}
\vskip 1cm

\vskip 1cm

\begin{proposition}\textbf{Harmonic extension of a function, on the graph of the Weierstrass function}\\

\noindent For any strictly positive integer~$m$, if~$u$ is a real-valued function defined on~$V_{m-1}$, its \textbf{harmonic extension}, denoted by~$ \tilde{u}$, is obtained as the extension of~$u$ to~$V_m$ which minimizes the energy:

$$  {\cal{E}}_{\Gamma_{{\cal W}_m}}(\tilde{u},\tilde{u})=\eta_{2-D_{\cal W}   }^{-2}\,N_b^{ 2\,( 2-D_{\cal W} )\,m}\,\displaystyle \sum_{X \underset{m }{\sim} Y} (\tilde{u}(X)-\tilde{u}(Y))^2 $$

\noindent The link between~$   {\cal{E}}_{\Gamma_{{\cal W}_m}}$ and~$  {\cal{E}}_{\Gamma_{{\cal W}_{m-1}}}$ is obtained through the introduction of two strictly positive constants~$r_m$ and~$r_{m+1}$ such that:

$$ \eta_{2-D_{\cal W}   }^{-2}\,N_b^{ 2\,(2- D_{\cal W} )\,m}\, r_{m }\, \displaystyle \sum_{X \underset{m  }{\sim} Y} (\tilde{u}(X)-\tilde{u}(Y))^2 =\eta_{2-D_{\cal W}   }^{-2}\,N_b^{ 2\,( 2-D_{\cal W} )\,(m-1)} r_{m-1}\,\displaystyle  \sum_{X \underset{m-1 }{\sim} Y} (u(X)-u(Y))^2$$

\noindent In particular:

$$  \eta_{2-D_{\cal W}   }^{-2}\,N_b^{ 2\,( 2-D_{\cal W} ) }\,r_{1 }\, \displaystyle \sum_{X \underset{1  }{\sim} Y} (\tilde{u}(X)-\tilde{u}(Y))^2 = \eta_{2-D_{\cal W}   }^{-2}\, r_{0}\,\displaystyle  \sum_{X \underset{0 }{\sim} Y} (u(X)-u(Y))^2$$

\noindent For the sake of simplicity, we will fix the value of the initial constant:~$r_0=1$. One has then:

$$ {\cal{E}}_{\Gamma_{{\cal W}_m}}(\tilde{u},\tilde{u})= \displaystyle \frac{1}{ r_{1 }}\,  {\cal{E}}_{\Gamma_{{\cal W}_0}}(\tilde{u},\tilde{u})$$

\noindent Let us set:

$$r = \displaystyle \frac{1}{r_{1 }} $$
%\noindent ce qui peut aussi s'écrire sous la forme :

%$$  \displaystyle \sum_{i=0}^{N_b-1} r_i^{-1}  {\cal{E}}_{\Gamma_{{\cal W}_{ m-1}}}(\tilde{u} \circ T_i,\tilde{u} \circ T_i )   =   {\cal{E}}_{\Gamma_{{\cal W}_{ m }}}(\tilde{u} %,\tilde{u}  )  $$

%\noindent où les~$r_i$,~\mbox{$0 \leq i \leq N_b-1$}, sont des constantes de normalisation.\\

\noindent and:

$$  {\cal{E}}_{m}(u)=\eta_{2-D_{\cal W}   }^{-2}\,r_m\, \sum_{X \underset{m }{\sim} Y} (\tilde{u}(X)-\tilde{u}(Y))^2 $$

\noindent Since the determination of the harmonic extension of a function appears to be a local problem, on the graph~$\Gamma_{{\cal W}_{ m-1}}$, which is linked to the graph~$\Gamma_{{\cal W}_{ m }}$ by a similar process as the one that links~$\Gamma_{{\cal W}_{  1}}$ to~$\Gamma_{{\cal W}_{ 0}}$, one deduces, for any strictly positive integer~$m$:

$$ {\cal{E}}_{\Gamma_{{\cal W}_m}}(\tilde{u},\tilde{u})= \displaystyle \frac{1}{ r_{1 }}\,  {\cal{E}}_{\Gamma_{{\cal W}_{m-1}}}(\tilde{u},\tilde{u})$$

\noindent By induction, one gets:

$$r_m=r_1^m\,r_0=r^{-m} $$

%\noindent La suite~$\left (\mathcal{E}_m(u)\right)_{m\in\N}$ est croissante.

\noindent If~$v$ is a real-valued function, defined on~$V_{m-1}$, of harmonic extension~$ \tilde{v}$, we will write:

$$  {\cal{E}}_{m}(u,v)=\eta_{2-D_{\cal W}   }^{-2}\,r^{-m}\, \sum_{X \underset{m }{\sim} Y} (\tilde{u}(X)-\tilde{u}(Y)) \, (\tilde{v}(X)-\tilde{v}(Y)) $$

\noindent For further precision on the construction and existence of harmonic extensions, we refer to~\emph{\cite{Sabot1987}}.
\end{proposition}

\vskip 1cm

\begin{pte}\textbf{Self-similar measure, for the domain delimited by the graph of the Weierstrass function}\\

\noindent Let us denote by~$  \mu_{\cal L}$ the Lebesgue measure on~$\R^2$. We set, for any~$i$ of~\mbox{$ \left \lbrace 0, \hdots, N_b-1 \right \rbrace$ }:
$$ \mu_i= \displaystyle \frac{  \mu_{\cal L}\left ( T_i\left ({\cal P}_0\right) \right) }{\mu_{\cal L} \left (  {\cal P}_0\right)}$$

\noindent The measure~$\mu $, such that:

$$ \mu= \displaystyle \sum_{i=0}^{N_b-1} \mu_i\,\mu\circ T_i^{-1}  $$

\noindent is self-similar, for the domain delimited by the graph of the Weierstrass function. We refer to~\cite{ClaireGB} for further details.

\end{pte}

\vskip 1cm

\begin{definition}\textbf{Laplacian of order~$m\,\in\,\N^\star$}\\

\noindent For any strictly positive integer~$m$, and any real-valued function~$u$, defined on the set~$V_m$ of the vertices of the graph~$\Gamma_{{\cal W}_m}$, we introduce the Laplacian of order~$m$,~$\Delta_m(u)$, by:

$$\Delta_m u(X) = \displaystyle\sum_{Y \in V_m,\,Y\underset{m}{\sim} X} \left (u(Y)-u(X)\right)  \quad \forall\, X\,\in\, V_m\setminus V_0 $$

\end{definition}

\vskip 1cm

\begin{definition}\textbf{Existence domain of the Laplacian, for a continuous function on the graph~$\Gamma_{\cal W}$ } (see \cite{Beurling1985})\\

\label{Lapl}
\noindent We will denote by~$\text{dom}\, \Delta$ the existence domain of the Laplacian, on the graph~$\Gamma_{\cal W}$, as the set of functions~$u$ of~$\text{dom}\, \mathcal{E}$such that there exists a continuous function on~$\Gamma_{\cal W}$, denoted~$\Delta \,u$, that we will call \textbf{Laplacian of~$u$}, such that :
$$\mathcal{E}(u,v)=-\displaystyle \int_{{\cal D} \left ( \Gamma_{\cal W} \right)} v\, \Delta u   \,d\mu \quad \text{for any } v \,\in \,\text{dom}_0\, \mathcal{E}$$
\end{definition}

\vskip 1cm

\begin{notation}

In the following, we will denote by~${\cal H}_0\subset \text{dom}\, \Delta$ the space of harmonic functions, i.e. the space of functions~$u \,\in\,\ \text{dom}\, \Delta$ such that:

$$\Delta\,u=0$$

\noindent Given a natural integer~$m$, we will denote by~${\cal S} \left ({\cal H}_0,V_m \right)$ the space, of dimension~$N_b^m$, of spline functions " of level~$m$", ~$u$, defined on~$\Gamma_{\cal W}$, continuous, such that, for any word~$\cal M$ of length~$m$,~\mbox{$u \circ T_{\cal M}$} is harmonic, i.e.:

$$\Delta_m \, \left ( u \circ T_{\cal M} \right)=0$$

\end{notation}

\vskip 1cm

\begin{pte}
Let~$m$ be a strictly positive integer,~$X \,\notin\,V_0$ a vertex of the graph~$\Gamma_{\cal W}$, and~\mbox{$\psi_X^{m}\,\in\,{\cal S} \left ({\cal H}_0,V_m \right)$} a spline  function such that:

$$\psi_X^{m}(Y)=\left \lbrace \begin{array}{ccc}\delta_{XY} & \forall& Y\,\in \,V_m \\
 0 & \forall& Y\,\notin \,V_m \end{array} \right. \quad,  \quad \text{where} \quad    \delta_{XY} =\left \lbrace \begin{array}{ccc}1& \text{if} & X=Y\\ 0& \text{else} &  \end{array} \right.$$

\noindent For any function~$u$ of~$\text{dom}\, \mathcal{E}$, such that its Laplacian exists:

$$  \Delta  u(X)= \displaystyle \lim_{m \to + \infty} \eta_{2-D_{\cal W}   }^{-2}\,r^{-m} \left (  \displaystyle\int_{{\cal D} \left ( \Gamma_{\cal W} \right)}  \psi_X^{m}\, d\mu  \right)^{-1} \,\Delta_m u(X)\,$$

\end{pte}

\vskip 1cm

\begin{notation}
\noindent We will denote by~$\text{dom}\,{\cal E}$ the subspace of continuous functions defined on~$\Gamma_{\cal W}$, such that:

$$\mathcal{E}(u)< + \infty$$

\end{notation}

\vskip 1cm

\begin{pte}\textbf{Spectrum of the Laplacian}(We refer to our work~\cite{ClaireGB})\\

\noindent Let us consider the eigenvalues~$\left (-{\Lambda_m}\right)_{m\in\N}$ of the sequence of graph Laplacians~$\left (\Delta_m  \right)_{m\in\N}$, built on the discrete sequence of graphs~$\left (\Gamma_{{ \cal W}_m}\right)_{m\in\N}$.\\

\noindent The spectral decimation method leads to the following recurrence relations between the eigenvalues of order~$m$ and~$m+1$:

$$ \Lambda_{m }
= \displaystyle \frac{\left \lbrace \left (\displaystyle \frac{-2+\Lambda_{m-1 }-\varepsilon \,\left (\left \lbrace \Lambda_{m-1 }-2 \right\rbrace ^2-4 \right)^{\frac{1}{2}}}{2}\right)
^{\frac{1}{N_b}}   +1 \right \rbrace^2}{\left (\displaystyle \frac{-2+\Lambda_{m-1 }-\varepsilon \,\left (\left \lbrace \Lambda_{m-1 }-2 \right\rbrace ^2-4\right)^{\frac{1}{2}}}{2}\right)
^{\frac{1}{N_b}}   }
$$

\noindent where~\mbox{$\varepsilon\,\in\,\left \lbrace -1, 1 \right \rbrace$}.

\end{pte}

\vskip 1cm

\section{Effective resistance metric, on the graph of the Weierstrass function}

\begin{pte}
\noindent The space~$\text{dom}\, \Delta$, modulo constant functions, is a Hilbert space, included in the space of continuous functions on the graph~$\Gamma_{\cal W}$, modulo constant functions.
\end{pte}

\vskip 1cm

%\begin{proof}A FAIRE.

%\end{proof}

%\vskip 1cm

\begin{definition}\textbf{Effective resistance metric, on the graph~$\Gamma_{\cal W}$}\\
\noindent Given a pair of points~$(X,Y)$ of the graph~\mbox{$ \Gamma_{\cal W} $}, we define, as in~\cite{Strichartz2003}, the \textbf{effective resistance metric between the points~$X$ and~$Y$}, by:

$$ R_{\Gamma_{\cal W}}(X,Y)= \left \lbrace \min_{\left \lbrace u \, |\, u(X)=0 , u(Y)=1\right \rbrace } \,{\cal E}(u )  \right \rbrace^{-1}$$

\noindent In an equivalent way,~\mbox{$ R_{\Gamma_{\cal W}}(X,Y)$} may be defined as the minimum value of the real numbers~$R $ such that, for any function~$u$ of~\mbox{$\text{dom} \,\Delta$}:

$$ \left |u (X)-u (Y)\right |^2 \leq R\,  {\cal E}(u)  $$

\end{definition}

\vskip 1cm

\begin{definition}\textbf{Metric, on the graph~$\Gamma_{\cal W}$}\\
\noindent Let us define, on the graph~$\Gamma_{\cal W}$, the distance~\mbox{$d_{\Gamma_{\cal W}}$} defined, for any pair of points~$(X,Y)$ of~\mbox{$ \Gamma_{\cal W} $}, by:

$$ d_{\Gamma_{\cal W}}(X,Y)= \left \lbrace \min_{\left \lbrace u \, |\, u(X)=0 , u(Y)=1\right \rbrace } \,{\cal E}(u,u)  \right \rbrace^{-1}$$
\end{definition}

\vskip 1cm

\begin{remark}

As it is explained in~\cite{StrichartzLivre2006}, one may note that the minimum

$$\min_{\left \lbrace u \, |\, u(X)=0 , u(Y)=1\right \rbrace } \,{\cal E}(u )   $$

\noindent is reached when the function~$u$ is harmonic on the complement set, in~\mbox{$\Gamma_{\cal W}$},
of the set~\mbox{$ \left \lbrace X \right \rbrace \cup  \left \lbrace Y \right \rbrace  $} (we recall that, by definition, a harmonic function~$u$ on~
\mbox{$\Gamma_{\cal W}$} minimizes the sequence of energies~\mbox{$\left ({\cal E}_{\Gamma_{ {\cal W}_m}}  (u,u )\right)_{m\in\N}$}.\\

\noindent In order to fully apprehend and understand the intrinsic meaning of these functions, one might reason by analogy with the unit interval~$[0,1]$. In this case, one will note that, given two points~$X$ and~$Y$ of~$[0,1]$ such that~\mbox{$X < Y$}, the function~$u$ is affine by pieces, taking the value zero on~$[0,X]$, and the value~1 on~$[Y,1]$ (see the illustration on the following figure):

$$\forall \,t \,\in\,[0,1] \, : \quad u(t)=\displaystyle \frac { t-X}{ Y-X}$$

 \begin{figure}[h!]
  \center{ \psfig{height=5cm,width=14cm,angle=0,file=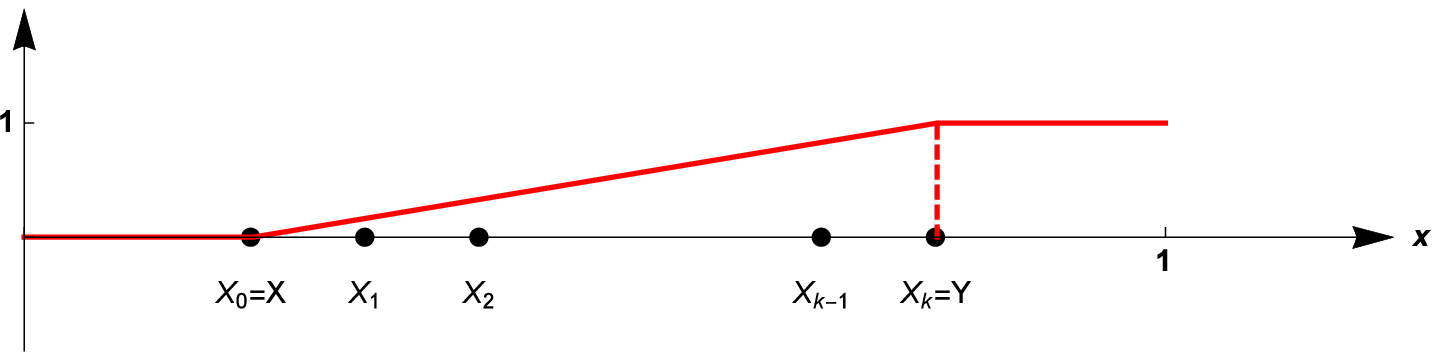}}
 \caption{The graph of the function~$u$ where the value~\mbox{$ \min_{\left \lbrace u \, |\, u(X)=0 , u(Y)=1\right \rbrace } \,{\cal E}(u )    $} is reached.}
 \end{figure}

\noindent Let us denote by~$m$ the natural integer such that:

$$X  \underset{m }{\sim}  Y $$

\noindent One may introduce, the, for any integer~$p$, the sequence of points~\mbox{$\left (X_{ j}\right)_{0 \leq j \leq 2^p}$} such that:

$$X_0= X \quad , \quad X_{2^p}=Y $$

\noindent and, for any integer~$j$ such that~\mbox{$0 < j < 2^p-1$}:

$$ X_{ j} \,\in\,V_{p+1} \quad , \quad X_{ j}  \underset{p+ 1 }{\sim}  X_{ j+1} $$

\noindent In the case of the unit interval, the normalization constant is:

$$r^{-1} = \displaystyle \frac{1}{2}$$

\noindent One has then:

$$\begin{array}{ccc} {\cal{E}} (u,u)&=&\displaystyle \lim_{p \to + \infty} {\cal{E}}_{p }(u,u)\\
&=&\displaystyle \lim_{m \to + \infty} r^{-p}\,{\cal E}_{\Gamma_{ {\cal W}_p}}  (u )\\
  &=& \displaystyle \lim_{p \to + \infty} \displaystyle \sum_{(X,Y) \,\in \, V_p^2,\,X  \underset{p }{\sim}  Y} r^{-p}\,\left (u_{\mid V_p}(X)-u_{\mid V_p}(Y)\right )^2 \\
  &=& \displaystyle \lim_{p \to + \infty} \displaystyle \sum_{(X,Y) \,\in \, V_p^2,\,X  \underset{p }{\sim}  Y} \displaystyle \frac{1}{2^p}\,\left (u_{\mid V_p}(X)-u_{\mid V_p}(Y)\right )^2 \\
  &=& \displaystyle \lim_{p \to + \infty} \displaystyle \sum_{j=0}^{ k-1}  \displaystyle \frac{1}{2^p}\,\left (u \left ( X+\displaystyle \frac{ j}{2^p}\right )-u \left ( X+\displaystyle \frac{ j+1}{2^p}\right )\right )^2 \\
  &=& \displaystyle \int_X^Y\displaystyle \frac { dt}{(Y-X)^2} \\
    &=&  \displaystyle \frac {1}{Y-X} \\
  \end{array}
  $$

  \noindent If~$d_{\R}$ denotes the usual Euclidean distance on~$\R$:

  $$\forall\, (X,Y)\,\in\R^2 \, : \quad d_{\R}(X,Y)=|Y-X|$$

  \noindent one has then:

$$\begin{array}{ccc} \displaystyle \min_{\left \lbrace u \, |\, u(X)=0 , u(Y)=1\right \rbrace } \,{\cal E}(u )  &=&     \displaystyle \frac {1}{d_{\R}(X,Y)}
  \end{array}
  $$

  \noindent Let us now consider, more generally, a fractal domain~$\cal F$, in an Euclidean space of dimension~\mbox{$d \,\in\,\N^\star$}, equipped with the distance~
  \mbox{$d_{\R^d}$}. If, one has, in advance, defined an energy on~$\cal F$, it is worth searching wether there exists a real number~$\beta$ such that:

$$\forall\, (X,Y) \,\in \, {\cal F}^2 \, : \quad  \left (\min_{\left \lbrace u \, |\, u(X)=0 , u(Y)=1\right \rbrace } \,{\cal E}(u )  \right)^{-1} \sim \left (d_{\R^d}(X,Y) \right)^\beta$$

 \noindent In the case of the Sierpi\`{n}ski gasket~$\cal SG$ (we refer to~\cite{StrichartzFunctionalAnalysis}), Robert~S.~Strichartz lays the emphasis upon the fact that, given~\mbox{$X  \underset{m }{\sim}  Y $}, one has:

 $$ \displaystyle \min_{\left \lbrace u \, |\, u(X)=0 , u(Y)=1\right \rbrace } \,{\cal E}(u ) \lesssim r_{\cal SG}^m=\left (\displaystyle \frac{3}{5} \right)^m$$

 \noindent This also corresponds thus to the order of the diameter of the~$m^{th}-$order cells.\\

\noindent Since the Sierpi\`{n}ski gasket~$\cal SG$ is obtained from the initial triangle of diameter~1 by means of three contractions, the respective ratios of which are equal to~$\displaystyle \frac{1}{2}$, one has simply to look the real number~$\beta_{\cal SG}$ such that:

 $$\left (\displaystyle \frac{1}{2} \right)^{m\,\beta_{\cal SG}}=\left (\displaystyle \frac{3}{5} \right)^m$$

 \noindent This leads to:

 $$ \beta_{\cal SG} = \displaystyle \frac{\ln \frac{5}{3}}{\ln 2}  $$

\end{remark}

\vskip 1cm

\newpage

\begin{definition}\textbf{Dimension of the graph~$\Gamma_{\cal W}$, in the effective resistance metric}\\

\noindent The \textbf{dimension of the graph~$\Gamma_{\cal W}$}, in the effective resistance metric, is the strictly positive number~\mbox{$d_{Gamma_{\cal W} }$} such that, given a strictly positive real number~$r$, and a point~$X\,\in\,{ \Gamma_{\cal W} } $, for the~\mbox{$X-$centered} ball of radius~$r$, denoted by~\mbox{${\cal B}_r(X)$}:

$$ \mu \left ({\cal B}_r(X) \right) =r^{d_{\Gamma_{\cal W}}}$$

\end{definition}

\vskip 1cm

\begin{proposition}
\noindent The dimension of the graph~$\Gamma_{\cal W} $, in the effective resistance metric, is given by:

\begin{enumerate}
\item[\emph{i}.] \textbf{First case:}~\mbox{$\lambda > \displaystyle \frac{1}{N_b}$}.\\

$$ d_{\Gamma_{\cal W}}=\displaystyle \frac{ \ln \frac{N_b}{  \lambda}}{(5-2\,D_{\cal W})\,\ln N_b} $$

%$$D_{\cal W}=2+ \displaystyle \frac{\ln \lambda}{\ln N_b}$$

\item[\emph{ii}.] \textbf{Second case:}~\mbox{$\lambda < \displaystyle \frac{1}{N_b}$}.\\

$$ d_{\Gamma_{\cal W}}=\displaystyle \frac{ 2}{ 5-2\,D_{\cal W} } $$

\end{enumerate}

\end{proposition}

\vskip 1cm

\begin{proof}

\begin{remark}

\noindent Once again, it is worth having a look at the case of the Sierpi\'{n}ski gasket. Robert~S.~Strichartz stars from the fact that the measure of~$m^{th}-$order cells is~$\displaystyle \frac{1}{3^m}$. Two consecutive points~$x$ and~$y$ are such that, for the effective resistance metric

$$d(x,y) \sim \left (\displaystyle \frac{3}{5 }\right)^m$$

\noindent For the self-similar measure~$\mu_{\cal SG}$, which affects the value~$\displaystyle \frac{1}{3^m}$ to each ~$m^{th}-$order cell, one has simply to look for the real number~$d_{\cal SG}$ such that:

$$\left (\displaystyle \frac{3}{5 }\right)^{m \,d_{\cal SG}}= \displaystyle \frac{1}{3^m}$$

\noindent which leads to:

$$d_{\cal SG} = \displaystyle \frac{\ln 3 }{\ln \frac{5}{3}}$$

\noindent One may then deduce from the above an estimate, for the effective resistance metric, of the measure of a~\mbox{$X-$centered} ball of radius~$r$, denoted by~\mbox{${\cal B}_r(X)$}:

$$\mu_{\cal SG} \left ({\cal B}_r(x) \right)=r^{d_{\cal SG}}$$

\end{remark}
\vskip 1cm

\noindent Let us now go back to the graph~$\Gamma_{\cal W}$.\\

  \noindent Given a natural integer~$m$, and two points~$X$ and~$Y$ such that~\mbox{$X  \underset{m }{\sim}  Y $}:

 $$ \displaystyle \min_{\left \lbrace u \, |\, u(X)=0 , u(Y)=1\right \rbrace } \,{\cal E}(u )\lesssim \eta_{2-D_{\cal W}   }^{-2} \,r^{-m}= \eta_{2-D_{\cal W}   }^{-2} \,
 N_b^{(5- 2\,D_{\cal W})\,m   }$$

\noindent For the detailed calculations which enable one to obtain the normalization constants, we refer to~\cite{DavidRiane}.

\noindent For the self-similar measure~$\tilde{\mu}$ introduced in the above, each~$m^{th}-$order cell, i.e. each simple polygon~${\cal P}_{m,j}$,~\mbox{$0 \leq j \leq N_b^m-1$}, with~$N_b$ sides and~$N_b$ vertices, has a measure of the order of:

 $$(N_b-1)\,  \displaystyle \frac{\eta^m}{N_b^m}  $$

\noindent The points~$X$ and~$Y$ such that~\mbox{$X  \underset{m }{\sim}  Y $} belong to a~$m^{th}-$order subcell, which is the intersection of a simple polygon~${\cal P}_{m,j}$,~\mbox{$0 \leq j \leq N_b^m-1$}, with the rectangle of which~$X$ and~$Y$ are two vertices, of width~\mbox{$\displaystyle \frac{\eta^m}{N_b^m} $}, and height~$\eta^m$. This subcell a has a measure, the order of which is thus:

 $$ \eta^m \,L_m = \displaystyle \frac{\eta^m}{(N_b-1)\,N_b^m}  $$
\begin{enumerate}
\item[\emph{i}.] \textbf{First case:}~\mbox{$\lambda > \displaystyle \frac{1}{N_b}$}.\\

\noindent One has simply to look for the real number~$d_{\Gamma_{\cal W}}$ such that:

$$\left (N_b^{5- 2\,D_{\cal W}   }\right)^{m\,  d_{\Gamma_{\cal W}}}  \quad \text{and} \quad  \displaystyle \frac{\lambda^m}{N_b^m}  $$
\noindent are of the same order, which yields:
$$ d_{\Gamma_{\cal W}}=\displaystyle \frac{ \ln \frac{N_b}{  \lambda}}{(5-2\,D_{\cal W})\,\ln N_b} $$

\item[\emph{ii}.] \textbf{Second case:}~\mbox{$\lambda < \displaystyle \frac{1}{N_b}$}.\\

\noindent One has simply to look for the real number~$d_{\Gamma_{\cal W}}$ such that:

$$\left ( N_b^{5- 2\,D_{\cal W}   } \right)^{m\,d_{\Gamma_{\cal W}}}      \quad \text{and} \quad  \displaystyle \frac{1}{N_b^{2m}}  $$
\noindent are of the same order, which yields:
$$ d_{\Gamma_{\cal W}}=\displaystyle \frac{ 2}{  5-2\,D_{\cal W}  } $$

\end{enumerate}
\end{proof}

% \noindent ~$\Gamma_{\cal W}$ est alors un espace homogène de dimension~\mbox{$d_{\Gamma_{\cal W}}$} pour la métrique de résistance et la mesure~$\tilde{\mu}$.

\section{Detailed study of the spectrum of the Laplacian}

%\subsection{Comparison with the classical graph Laplacian matrix}

%$${\cal L}_{V_0}=D_{V_0}-A_{V_0}=
% \left ( \begin{array}{ccc} 1&0&0 \\0&2&0 \\0&0&1 \end{array} \right)-
% \left ( \begin{array}{ccc} 0&1&0 \\1&0&1 \\0&1&0 \end{array} \right)
 % =\left ( \begin{array}{ccc} 1&1&0 \\1&2&1 \\0&1&1 \end{array} \right)$$

 % $${\cal L}_{V_1}=D_{V_1}-A_{V_1}=
% \left ( \begin{array}{ccccccc} 1&0&0 &0&0&0&0\\0&2&0 &0&0&0&0\\0&0&2&0&0&0&0\\0&0 &0&2&0&0&0 \\0&0&0&0&2&0&0 \\0&0&0&0&0&2&0 \\0&0&0&0&0&1&2  \end{array} \right)-
 % \left ( \begin{array}{ccccccc} 0&1&0&0&0&0&0 \\1&0&1&0&0&0&0 \\0&1&0&1&0&0&0 \\0&0&1&0&1&0&0\\0&0&0&1&0 &1&0  \\0&0&0&0&1 &0&1   \\0&0&0&0&0 &1&0  \end{array} \right)
%  = \left ( \begin{array}{ccccccc} 1&1&0&0&0&0&0 \\1& 2&1&0&0&0&0 \\0&1&2&1&0&0&0 \\0&0&1&2&1&0&0\\0&0&0&1&2 &1&0  \\0&0&0&0&1 &2&1   \\0&0&0&0&0 &1&1  \end{array} \right)
%  $$

 \noindent As exposed by R.~S.~Strichartz in~\cite{StrichartzLivre2006}, one may bear in mind that the eigenvalues can be grouped into two categories:

 \begin{enumerate}

 \item[\emph{i}.]  initial eigenvalues, which a priori belong to the set of forbidden values (as for instance~\mbox{$\Lambda=2$}) ;

 \item[\emph{ii}.] continued eigenvalues, obtained by means of spectral decimation.

 \end{enumerate}

 \vskip 1cm

 \noindent We present, in the sequel, a detailed study of the spectrum of~$\Delta$, in the case where~$N_b=3$, which can be easily extended to higher values of the integer~$N_b$.
 \vskip 1cm

 \subsection{Eigenvalues and eigenvectors of~$\Delta_1$}

\noindent Let us recall that the vertices of the graph~$\Gamma_{{\cal W}_1}$ are:

$$P_0 \quad , \quad T_0\left (P_1 \right)\quad , \quad T_0\left (P_2 \right) \quad , \quad T_1\left (P_0 \right)$$

$$P_1\quad , \quad T_1\left (P_2 \right)\quad , \quad T_2\left (P_0 \right) \quad , \quad T_2\left (P_1 \right)\quad , \quad P_2$$

\noindent One may note that:

$$\text{Card}\, \left ( V_1\setminus V_0 \right)=4$$

\noindent Let us denote by~$u$ an eigenfunction, for the eigenvalue~$-\Lambda$. For the sake of simplicity, we set:

$$u\left (T_0\left (P_1 \right)\right)=a\,\in\,\R \quad, \quad u\left (T_0\left (P_2 \right)\right)=b\,\in\,\R \quad, \quad
u\left (T_2\left (P_0 \right)\right)=c\,\in\,\R \quad, \quad  u\left (T_2\left (P_1 \right)\right)=d\,\in\,\R $$

 \begin{figure}[h!]
 \center{\psfig{height=8cm,width=14cm,angle=0,file=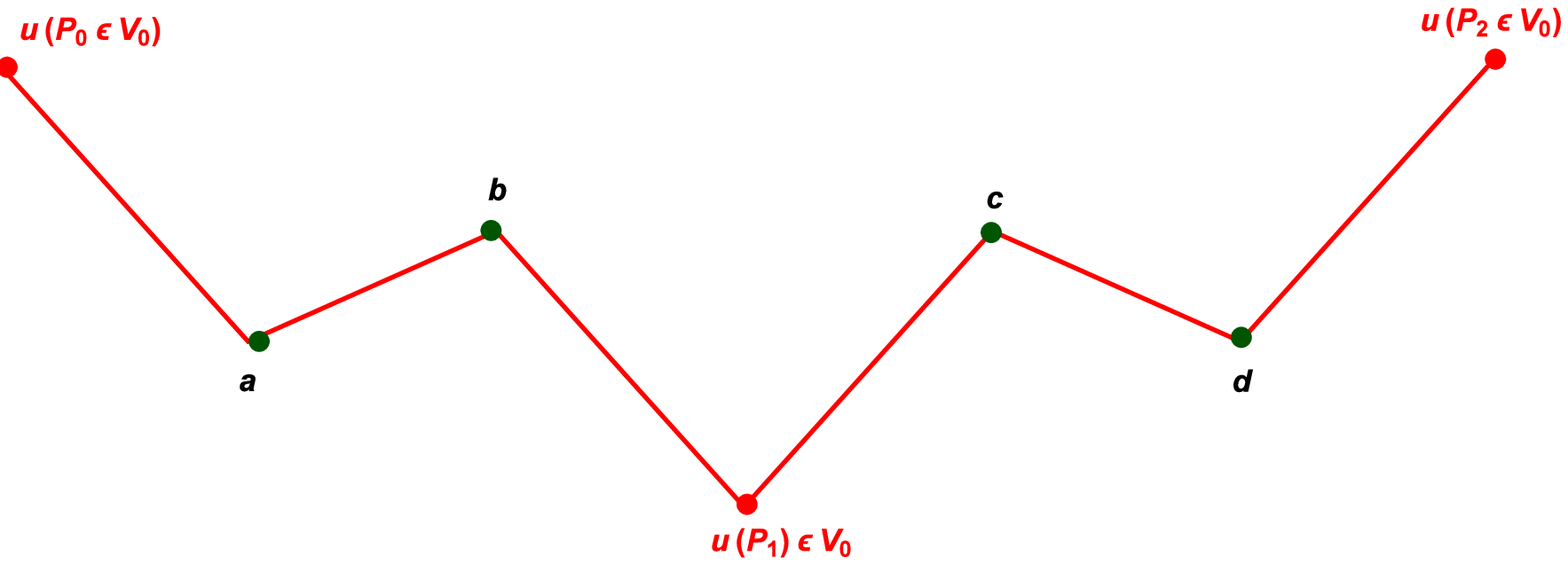}}
\caption{Successive values of an eigenfunction on~$V_1$, in the case where~$N_b=3$.}
 \end{figure}

\noindent One has then:

$$\left \lbrace \begin{array}{ccc}
u(P_0)-a+b-a &=& -\Lambda\,a \\
a-b+u(P_1)-b&=& -\Lambda\,b \\
u(P_1)-c+d-c&=& -\Lambda\,c \\
c-d+u(P_2)-d&=& -\Lambda\,d \end{array} \right.$$

\noindent One may note that the only "Dirichlet eigenvalues", i.e. the ones related to the Dirichlet problem:

$$u_{|V_0}=0 \quad \text{i.e.} \quad u(P_0)=u(P_1)=u(P_2)=0$$

\noindent are obtained for:

$$\left \lbrace \begin{array}{ccc}
   b  &=& -(\Lambda-2)\,a \\
a   &=& -(\Lambda-2)\,b \\
 d &=& -(\Lambda-2)\,c \\
c   &=& -(\Lambda-2)\,d \end{array} \right.$$

\noindent i.e.:

$$\left \lbrace \begin{array}{ccc}
   b  &=&  (\Lambda-2)^2\,b \\
a   &=&  (\Lambda-2)^2\,a \\
 d &=&  (\Lambda-2)^2\,d \\
c   &=&  (\Lambda-2)^2\,c \end{array} \right.$$

\noindent The forbidden eigenvalue~$\Lambda=2$ cannot thus be a Dirichlet one.\\

\vskip 1cm
\noindent Let us consider the case where:

$$ (\Lambda-2)^2=1$$

\noindent i.e.

$$  \Lambda =1 \quad \text{or} \quad \Lambda =3   $$

\noindent The value~$\Lambda= 1$ leads to:

$$a=b\quad , \quad c=d$$

\noindent which yields a two-dimensional eigenspace. The multiplicity of the eigenvalue~$\Lambda =1$ is~2.\\\\

\noindent For the eigenvalue~$   \Lambda =3$:

$$a =-b\quad, \quad c=-d$$

%\noindent the symmetry of the problem requires:

%$$a=d$$

\noindent The eigenspace, for the eigenvalue~$3$, has dimension~2. The multiplicity of the eigenvalue~$\Lambda =3$ is~2.\\

 \noindent Since the cardinal of~$V_1\setminus V_0 $ is:

$$ {\cal N}^{\cal S}_{ 1}-N_b=  2\, N_b  -2 =4$$

\noindent one may note that we have the complete spectrum.

\vskip 1cm

 \subsection{Eigenvalues of~$\Delta_2$}

 \noindent Let us now look at the spectrum of~$\Delta_2$. For the sake of simplicity, we will denote by~$a$,~$b$,~$c$,~$d$,~$e$,~$f$,~$g$,~$h$, the successive values of an eigenfunction at the~\mbox{$N_b^2-1$} points between~$P_0$ and~$P_1$,  and by~$a'$,~$b'$,~$c'$,~$d'$,~$e'$,~$f'$,~$g'$,~$h'$, the successive values of an eigenfunction at the~\mbox{$N_b^2-1$} points between~$P_1$ and~$P_2$, as it appears on the following figure.

 \begin{figure}[h!]
 \center{\psfig{height=8cm,width=14cm,angle=0,file=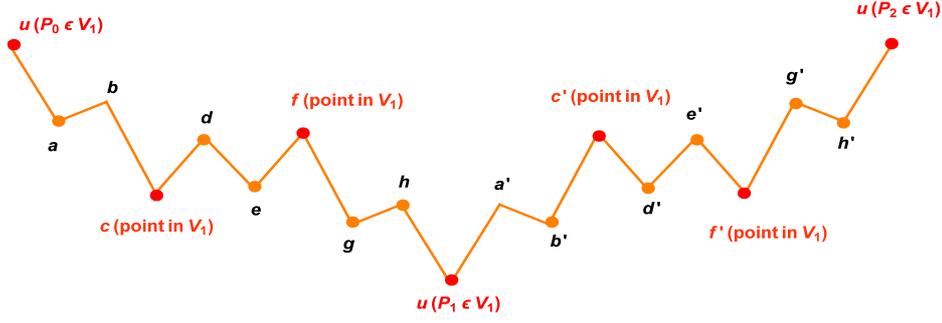}}
\caption{Successive values of an eigenfunction on~$V_2$, in the case where~$N_b=3$.}
 \end{figure}

\noindent One has then:

$$\left \lbrace \begin{array}{ccc}
 (2-\Lambda )\,a  &=& -u(P_0)-b\\
 (2-\Lambda )\,b  &=& -c-a\\
 (2-\Lambda )\,c  &=& -b-d\\
 (2-\Lambda )\,d  &=& -c-e\\
 (2-\Lambda )\,e  &=& -d-f\\
 (2-\Lambda )\,f  &=& -e-g\\
 (2-\Lambda )\,g  &=& -f-h\\
 (2-\Lambda )\,h  &=& -g-u(P_1)\\
\end{array} \right.$$

\noindent and:

$$ (2-\Lambda )\,u(P_1)   =  -h-a'  $$

\noindent and:

$$\left \lbrace \begin{array}{ccc}
 (2-\Lambda )\,a'  &=& -u(P_1)-b'\\
 (2-\Lambda )\,b'  &=& -c'-a'\\
  (2-\Lambda )\,c'  &=& -b'-d'\\
 (2-\Lambda )\,d'  &=& -c'-e'\\
 (2-\Lambda )\,e'  &=& -d'-f'\\
  (2-\Lambda )\,f'  &=& -e'-g'\\
 (2-\Lambda )\,g'  &=& -f'-h'\\
 (2-\Lambda )\,h'  &=& -g'-u(P_2)\\
\end{array} \right.$$

\noindent One may note that the only Dirichlet eigenvalues, in the case where:

$$u_{|V_0  }=0 \quad \text{i.e.} \quad u(P_0)=u(P_1)=u(P_2)= 0$$

\noindent are obtained for:

$$\left \lbrace \begin{array}{ccc}
 (2-\Lambda )\,a  &=&  -b\\
 (2-\Lambda )\,b  &=&  -a\\
  (2-\Lambda )\,c  &=& -b-d\\
 (2-\Lambda )\,d  &=&  -e\\
 (2-\Lambda )\,e  &=& -d \\
   (2-\Lambda )\,f   &=& -e-g\\
 (2-\Lambda )\,g  &=&  -h\\
 (2-\Lambda )\,h  &=& -g \\
\end{array} \right.  \quad \text{and} \quad
0   =  -h-a'   \quad \text{and} \quad \left \lbrace \begin{array}{ccc}
 (2-\Lambda )\,a'  &=&  -b'\\
 (2-\Lambda )\,b'  &=&  -a'\\
  (2-\Lambda )\,c'  &=& -b'-d'\\
 (2-\Lambda )\,d'  &=&  -e'\\
 (2-\Lambda )\,e'  &=& -d' \\
   (2-\Lambda )\,f'   &=& -e'-g'\\
 (2-\Lambda )\,g'  &=&  -h'\\
 (2-\Lambda )\,h'  &=& -g' \\
\end{array} \right.$$

\noindent i.e.:

$$\left \lbrace \begin{array}{ccc}
 (2-\Lambda )\,a  &=&  -b\\
 (2-\Lambda )^2\,b  &=&   b\\
 \left \lbrace 1-(2-\Lambda )^2\right \rbrace \,a  &=& 0 \\
   (2-\Lambda )\,c  &=& -b-d\\
 (2-\Lambda )\,d  &=&  -e\\
 (2-\Lambda )^2\,e  &=& e \\
   (2-\Lambda )\,f  &=& -e-g\\
 (2-\Lambda )^2\,h  &=&   h\\
 (2-\Lambda )\,h  &=& -g \\
\end{array} \right. \quad \text{and} \quad
0   =  -h-a'
 \quad \text{and} \quad \left \lbrace \begin{array}{ccc}
 (2-\Lambda )\,a'  &=&  -b'\\
 (2-\Lambda )^2\,b'  &=&   b'\\
 \left \lbrace 1-(2-\Lambda )^2\right \rbrace \,a'  &=& 0 \\
   (2-\Lambda )\,c'  &=& -b'-d'\\
 (2-\Lambda )\,d'  &=&  -e'\\
 (2-\Lambda )^2\,e'  &=& e' \\
   (2-\Lambda )\,f'   &=& -e'-g'\\
 (2-\Lambda )^2\,h'  &=&   h'\\
 (2-\Lambda )\,h'  &=& -g' \\
\end{array} \right.$$

\noindent The forbidden eigenvalue~$\Lambda=2$ is not therefore a Dirichlet one.\\

\vskip 1cm
\noindent Let us consider the case where:

$$ (\Lambda-2)^2=1$$

\noindent i.e.

$$  \Lambda =3 \quad \text{or} \quad \Lambda =1   $$

\noindent For~$\Lambda =1$, one has:

$$\left \lbrace \begin{array}{ccc}
  a  &=&  -b\\
    c  &=& -b-d\\
  d  &=&  -e\\
   f  &=& -e-g\\
 h  &=& -g \\
\end{array} \right. \quad \text{and} \quad
0   =  -h-a'
 \quad \text{and} \quad \left \lbrace \begin{array}{ccc}
 a'  &=&  -b'\\
   c'  &=& -b'-d'\\
 d'  &=&  -e'\\
    f'   &=& -e'-g'\\
  h'  &=& -g' \\
\end{array} \right.$$

\noindent The eigenspace, for~$\Lambda =1$, has thus dimension~5. The multiplicity of the eigenvalue~$\Lambda =1$ is~5.\\

\noindent For~$\Lambda =3$:

$$\left \lbrace \begin{array}{ccc}
 \left \lbrace 1-(2-\Lambda )^2\right \rbrace \,a  &=& 0 \\
    c  &=&  b+d\\
  d  &=&   e\\
   f  &=&  e+g\\
  h  &=&  g \\
\end{array} \right. \quad \text{and} \quad
0   =  -h-a'
 \quad \text{and} \quad \left \lbrace \begin{array}{ccc}
 a'  &=&   b'\\
    c'  &=&  b'+d'\\
  d'  &=&   e'\\
   f'   &=&  e'+g'\\
  h'  &=&  g' \\
\end{array} \right.$$

\noindent The eigenspace, for~$\Lambda =3$, has thus dimension~5. The multiplicity of the eigenvalue~$\Lambda =3$ is~5.\\

\noindent Let us now look at the continued eigenvalues, i.e. the ones obtained from the eigenvalues~\mbox{$\Lambda_1=1$} and~\mbox{$\Lambda_1=3$} by means of spectral decimation:

$$ \Lambda_{2 } = \phi^{-1} \left (\left (\phi\left (\Lambda_{1 } \right ) \right )^{\frac{1}{N_b}} \right)
=\displaystyle \frac{\left \lbrace \left (\phi\left (\Lambda_{1 } \right ) \right )^{\frac{1}{N_b}}  +1 \right \rbrace^2}{\left (\phi\left (\Lambda_{1 } \right )
\right )^{\frac{1}{N_b}}   }
= \displaystyle \frac{\left \lbrace \left (\displaystyle \frac{-2+\Lambda_{1 }-\varepsilon \,\sqrt{\left \lbrace \Lambda_{1 }-2 \right\rbrace ^2-4}}{2}\right)
^{\frac{1}{N_b}}   +1 \right \rbrace^2}{\left (\displaystyle \frac{-2+\Lambda_{1 }-\varepsilon \,\sqrt{\left \lbrace \Lambda_{1 }-2 \right\rbrace ^2-4}}{2}\right)
^{\frac{1}{N_b}}   }
$$

\noindent where~\mbox{$\varepsilon\,\in\,\left \lbrace -1, 1 \right \rbrace$}, for the values:

$$\Lambda_1 \,\in\, \left \lbrace 1,3 \right \rbrace $$

 \noindent As in~\cite{StrichartzLivre2006}, let us get rid, temporarily, of the Dirichlet conditions. We have thus:

$$\left \lbrace \begin{array}{ccc}
u(P_0) +b  &=& -(\Lambda-2)\,a \\
a +c &=& - (\Lambda-2)\,b \\
b+ d &=& -(\Lambda-2)\,c \\
e +f &=& -\Lambda\,d \\
e +g &=&-(\Lambda-2)\,f \\
f +h &=&-(\Lambda-2)\,g \\
g +u(P_1)&=& -(\Lambda-2)\,h \\
\end{array} \right. \quad \text{and} \quad (2-\Lambda )\,u(P_1)   =  -h-a'
  \quad \text{and} \quad  \left \lbrace \begin{array}{ccc}
u(P_1) +b'  &=& -(\Lambda-2)\,a' \\
a' +c' &=& - (\Lambda-2)\,b' \\
b'+ d' &=& -(\Lambda-2)\,c' \\
e' +f' &=& -\Lambda\,d' \\
e' +g' &=&-(\Lambda-2)\,f' \\
f' +h' &=&-(\Lambda-2)\,g' \\
g' +u(P_2)&=& -(\Lambda-2)\,h' \\
\end{array} \right.$$

\noindent For the initial eigenvalue~\mbox{$\Lambda_1=1$}, it is worth noticing that the restriction of the associated eigenvalue~$\Lambda$ to~$V_1\setminus V_0$ must satisfy the eigensystem associated to the eigenvalue~\mbox{$\Lambda_1=1$}, i.e.:

$$\left \lbrace \begin{array}{ccc}
u(P_0)+f&=& -(\Lambda_1-2)\,c \\
u(P_1)+c &=& -(\Lambda_1-2)\,f \\
\end{array} \right.  \quad \text{and} \quad  -( \Lambda-2 )\,u(P_1)   =   f+c'    \quad \text{and} \quad
\left \lbrace \begin{array}{ccc}
u(P_1)+f'&=& -(\Lambda_1-2)\,c' \\
u(P_2)+c' &=& -(\Lambda_1-2)\,f' \\
\end{array} \right.$$

\noindent or:
$$\left \lbrace \begin{array}{ccc}
u(P_0)+f&=& c \\
u(P_1)+c &=& f \\
\end{array} \right.   \quad \text{and} \quad  -( \Lambda-2 )\,u(P_1)   =   f+c'    \quad \text{and} \quad  \left \lbrace \begin{array}{ccc}
u(P_1)+f'&=& c' \\
u(P_2)+c' &=& f' \\
\end{array} \right.$$

\noindent i.e.:
$$
u(P_0)+u(P_1)=  0
\quad \text{and} \quad u(P_1)+u(P_2)=  0 $$

\noindent For~$u(P_0)=u(P_1)=u(P_2)=0$, it works, and the Dirichlet conditions appear to be satisfied. One has then:

\footnotesize

$$\left \lbrace \begin{array}{ccc}
 b  &=& -(\Lambda-2)\,a \\
 c &=& \left \lbrace- 1+ (\Lambda-2)^2\right \rbrace \,a \\
  d &=& (\Lambda-2)\,\left \lbrace 1 - \left \lbrace 1- (\Lambda-2)^2\right \rbrace \right \rbrace \,a \\
e +f &=& -\Lambda\,(\Lambda-2)\,\left \lbrace 1 - \left \lbrace 1- (\Lambda-2)^2\right \rbrace \right \rbrace \,a \\
e   &=&(\Lambda-2)\,\,\left \lbrace 1- \left \lbrace -1+(\Lambda-2)^2 \right \rbrace \right \rbrace \,h \\
f   &=&  \left \lbrace -1+(\Lambda-2)^2 \right \rbrace  \,h \\
g  &=& -(\Lambda-2)\,h \\
\end{array} \right.  \quad , \quad  c'  =  - f
 \quad, \quad  \left \lbrace \begin{array}{ccc}
 b'  &=& -(\Lambda-2)\,a' \\
 c' &=& \left \lbrace- 1+ (\Lambda-2)^2\right \rbrace \,a' \\
  d' &=& (\Lambda-2)\,\left \lbrace 1 - \left \lbrace 1- (\Lambda-2)^2\right \rbrace \right \rbrace \,a' \\
e' +f' &=& -\Lambda\,(\Lambda-2)\,\left \lbrace 1 - \left \lbrace 1- (\Lambda-2)^2\right \rbrace \right \rbrace \,a' \\
e'   &=&(\Lambda-2)\,\,\left \lbrace 1- \left \lbrace -1+(\Lambda-2)^2 \right \rbrace \right \rbrace \,h' \\
f'   &=&  \left \lbrace -1+(\Lambda-2)^2 \right \rbrace  \,h' \\
g'  &=& -(\Lambda-2)\,h' \\
\end{array} \right.$$

\normalsize
\noindent We obtain thus an eigenspace, the dimension of which is~$1$.\\

\noindent For the eigenvalue~$\Lambda_1=1$, the spectral decimation spectral leads to:

$$\left ( \displaystyle \frac{2- \Lambda_{ 2}+\varepsilon_2\, \rho\left (\omega_{ 2 } \right)^2 \, e^{  i \,\theta  \omega_{ 2 }} }{2}\right)^{N_b}
=\displaystyle \frac{1+\varepsilon_{ 1}\, \sqrt{3} \, e^{ \frac{i\,\pi}{2}} }{2}$$

\noindent which leads to the quadruple eigenvalue:

$$ \Lambda_{ 2}= 2+\cos \displaystyle \frac{ \pi}{9} +\sqrt{3}\, \sin\displaystyle \frac{ \pi}{9} $$

\noindent For the eigenvalue~$\Lambda_1=3$, the spectral decimation leads to the quadruple eigenvalue:

$$ \Lambda_{ 2}=  2\, \left \lbrace 1+ \cos \displaystyle \frac{ \pi}{9}   \right \rbrace $$

\noindent Since the cardinal of~$  V_2\setminus V_0 $ is:

 $$  {\cal N}^{\cal S}_{2}-{\cal N}^{\cal S}_{0}= 16$$

\noindent one may note that we have the complete spectrum.

\vskip 1cm

 \subsection{Eigenvalues of~$\Delta_3$}

 \noindent As previously, one can easily check that the forbidden eigenvalue~$\Lambda=2$ is not therefore a Dirichlet one.\\

\noindent One can also check that~$ \Lambda_3 =1$ and~$ \Lambda_3 =3$ are eigenvalues of~$\Delta_3$, both with multiplicity~8.\\

 \noindent From:

 $$ \Lambda_{  \hookrightarrow  2,}=  2\, \left \lbrace 1+ \cos \displaystyle \frac{ \pi}{9}   \right \rbrace $$

 \noindent the spectral decimation leads then to the quadruple eigenvalue:

 $$ \Lambda_{\hookrightarrow 3}= 4\, \cos^2   \displaystyle \frac{ \pi}{27} $$

 \noindent From:

 $$ \Lambda_{\hookrightarrow  2}= 2+\cos \displaystyle \frac{ \pi}{9} +\sqrt{3}\, \sin\displaystyle \frac{ \pi}{9}$$

 \noindent the spectral decimation leads then to the quadruple eigenvalue:

 $$ \Lambda_{ \hookrightarrow 3}= 4\, \cos^2   \displaystyle \frac{ \pi}{54} $$

\vskip 1cm

 \subsection{Eigenvalues of~$\Delta_m$,~$m \,\in\,\N$,~$m \geq 4$}

 \noindent As previously, one can easily check that the forbidden eigenvalue~$\Lambda=2$ is not therefore a Dirichlet one.\\

\noindent One can also check that~$ \Lambda_m =1$ and~$ \Lambda_m =3$ are eigenvalues of~$\Delta_m$, both with multiplicity~2.\\

 \noindent By induction, one may note that, due to the spectral decimation, the initial eigenvalue~$\Lambda_1=1$ gives birth, at this~$m^{th}$ step, to an eigenvalue~$\Lambda_m$, of multiplicity~$2^{m-1}$. In the same way, the initial eigenvalue~$\Lambda_1=3$ gives birth, at this~$m^{th}$ step, to an eigenvalue~$\Lambda_m$, of multiplicity~$2^{m-1}$.\\

 \noindent Results are summarized in the following array:

$$
\begin{array}{|c|c|c|c|c|r|}
\hline \\
\text{Initial eigenvalue}\, \Lambda_1 &\text{continued eigenvalue}\, \Lambda_2 &\text{continued eigenvalue $\Lambda_{\hookrightarrow 3}$}   &\text{continued eigenvalue  $\Lambda_{\hookrightarrow 4}$}\\
\hline \\
 1   & 2+\cos \displaystyle \frac{ \pi}{9} +\sqrt{3}\, \sin\displaystyle \frac{ \pi}{9}     &4\, \cos^2   \displaystyle \frac{ \pi}{27}   &4\, \cos^2   \displaystyle \frac{ \pi}{81}\\
\hline \\
3   &   2\, \left \lbrace 1+ \cos \displaystyle \frac{ \pi}{9}   \right \rbrace &4\, \cos^2   \displaystyle \frac{ \pi}{54} &
 2\, \left \lbrace 1+ \cos \displaystyle \frac{ \pi}{81}   \right \rbrace \\
\hline
\end{array}
$$

\vskip 1cm

\begin{pte}
\noindent Let us introduce:

$$\Lambda =  \displaystyle \lim_{m \to + \infty} \eta_{2-D_{\cal W}   }^{-2}\,  N_b^{ (5-2\,   D_{\cal W} )\,m } $$

\noindent One may note that, due to the definition of the Laplacian~$\Delta$, the limit exists.

\end{pte}

\vskip 1cm

 \subsection{Eigenvalue counting function}

 \begin{definition}\textbf{Eigenvalue counting function}\\

 \noindent Let us introduce the eigenvalue counting function, related to~$\Gamma_{\cal W} \setminus V_0$, such that, for any positive number~$x$:

  $${\cal N}^{\Gamma_{\cal W}\setminus V_0}(x) = \text{Card}\, \left \lbrace \Lambda \, \text{Dirichlet eigenvalue of~$-\Delta$} \, : \quad \Lambda \leq  x \right \rbrace $$

 \end{definition}

\vskip 1cm

 \begin{pte}

\noindent Given a strictly positive integer, the cardinal of~$V_m \setminus V_{0}$ is:

$$  {\cal N}^{\cal S}_{m } -{\cal N}^{\cal S}_{0 }
=2\, N_b^m -2 $$

%\noindent Let us denote by~$n_{m,1}$ and~$n_{m,3}$ the respective multiplicities of the initial eigenvalues~$\Lambda_m=1$ and~$\Lambda_m=3$. Previous results yield:

%$$ n_{m,1}=n_{m,3}=2$$

%\noindent Thus, the cardinal of the set of continued eigenvalues of~$-\Delta_m$ is:

%As the cardinal of~$V_m \setminus V_{m-1}$ is:

%$$  {\cal N}^{\cal S}_{m } -{\cal N}^{\cal S}_{m-1 }
%=2\,\left ( N_b^{m-1}- N_b^{m-2} \right) $$

 %\noindent Eigenvalue check:

% $$2 \times 2^{m-1}  +n_{m,1}+n_{m,3}=  {\cal N}^{\cal S}_{m } -{\cal N}^{\cal S}_{m-1 }
%=2\,\left ( N_b^{m-1}- N_b^{m-2} \right) $$

% $$  2^{m-1 }  +2\, n_{m,1}
%= \left ( N_b^{m-1}- N_b^{m-2} \right) $$

% $$     n_{m,1} =
%= \left ( N_b^{m-1}- N_b^{m-2} \right) -2^{m -1}  $$

\noindent Let us denote by~$\Lambda_m^s$ the largest eigenvalue, which is such that:

$$\eta_{2-D_{\cal W}   }\,N_b^m \times 3 \leq \Lambda_m^s \leq \eta_{2-D_{\cal W}   }\,N_b^m \times 4$$

\noindent  This leads to:

$$ {\cal N}^{\Gamma_{\cal W}} (C\, N_b^m   )= 2\, N_b^{m }- 2 \quad , \quad 3 \leq C \leq 4$$

 \noindent If one looks for an asymptotic growth rate of the form

$$ {\cal N}^{\Gamma_{\cal W}} (x) \sim x^\alpha$$

\noindent one obtains:

$$\alpha=1$$

\noindent By following~\cite{StrichartzLivre2006}, one may note that the ratio

$$ \displaystyle \frac{{\cal N}^{\Gamma_{\cal W}} (x)}{x}$$

\noindent is bounded above and away from zero, and admits a limit along any sequence of the form~$C \,N_b^m$,~\mbox{$C >0$},~\mbox{$m \,\in\,\N^\star$}. This enables one to deduce the existence of a periodic function~$g$, the period of which is equal to~$\ln N_b$, discontinuous at the value~$3$, such that:

$$\displaystyle \lim_{x \to + \infty} \left \lbrace \displaystyle \frac{{\cal N}^{\Gamma_{\cal W}} (x)}{x}-g (\ln x ) \right \rbrace =0$$

 \end{pte}

$$
\begin{array}{|c|c|c| }
\hline \\
\text{Level}  & \text{Cardinal of the Dirichlet spectrum, in the case where~$N_b=3$}     \\
\hline \\
 m   &  2\,  N_b^{m }- 2   \\
\hline \\
1   &  4 \\
\hline \\
2   &  16 \\
\hline \\
3   &   52 \\
\hline \\
4   &   160 \\
\hline
\end{array}
$$

 \begin{figure}[h!]
 \center{\psfig{height=8cm,width=14cm,angle=0,file=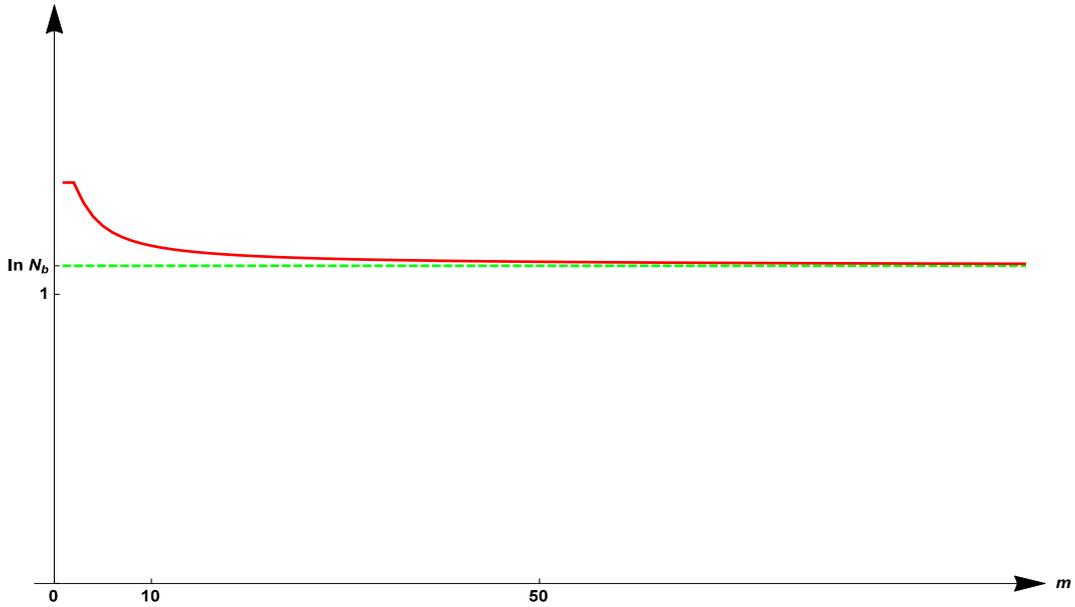}}\\
\caption{The graph of the function~$m \mapsto \displaystyle \frac{ln {{\cal N}^{\Gamma_{\cal W}}(C\,N_b^m)}}{m}$ in the case where~$\lambda= \displaystyle \frac{1}{2}$, and~$N_b=3$.}
 \end{figure}

 \vskip 1cm

 \begin{remark}

 \noindent Existing results of J.~Kigami and M.~Lapidus~\cite{KigamiLapidus1993}, and also of R.~S.~Strichartz~\cite{StrichartzLivre2006}, yield:

$$ {\cal N}^{\Gamma_{\cal W}} (x)= G(x)\, x^{\alpha_{\Gamma_{\cal W}}}+ {\cal O}(1)$$

\noindent with:

$$\alpha_{\Gamma_{\cal W}}=\displaystyle \frac{d_{\Gamma_{\cal W}}}{d_{\Gamma_{\cal W}}+1}
= \displaystyle \frac{\ln \frac  {N_b}{\lambda}} { \ln \frac  {N_b}{\lambda} +(5-2\,D_{\cal W})\,\ln N_b }$$

\noindent where:

$$ d_{\Gamma_{\cal W}}=\displaystyle \frac{ \ln \frac{N_b}{  \lambda}}{(5-2\,D_{\cal W})\,\ln N_b} $$

\noindent is the dimension of the graph~$\Gamma_{\cal W}$ for the resistance metric.\\

 \end{remark}
%\noindent This infirms the result of Strichartz ...\\

%$$C\, N_b^{m \,\alpha} \times 3^\alpha =  N_b^{m-1}- N_b^{m-2}$$

%$$\alpha=1$$

%$$C\,3^\alpha=N_b^{ -1}$$

\vskip 2cm

\centerline{\textbf{Thanks}}

\vskip 1cm
The author would like to thank R. Str., who suggested for our previous work, the introduction of specific energies to fully take into account the very specific geometry of the problem.


\vskip 2cm



\begin{thebibliography}{99}






\bibitem{Kigami1989}    J. Kigami, A harmonic calculus on the Sierpi\' nski spaces, Japan J. Appl. Math., \textbf{8} (1989), pages 259-290.

\bibitem{Kigami1993}  J. Kigami, Harmonic calculus on p.c.f. self-similar sets, Trans. Amer. Math. Soc., \textbf{335}(1993), pages 721-755.


\bibitem{Strichartz1999} R. S. Strichartz, Analysis on fractals, Notices of the AMS, \textbf{46}(8), 1999, pages 1199-1208.



\bibitem{Strichartz2} J. Kigami, R. S. Strichartz, K. C.  Walker, Constructing  a  Laplacian  on the  Diamond  Fractal, A. K.~ Peters,
 Ltd, Experimental  Mathematics, \textbf{10}(3), pages 437-448.


\bibitem{ClaireGB} Cl. David, Laplacian, on the graph of the Weierstrass function, arXiv:1703.03371v1.

\bibitem{Kaplan1984} J. Kaplan, J. Mallet-Paret and J. Yorke, The Lyapunov dimension of a nowhere differentiable attracting torus, Ergodic Theory Dynam. Systems
\textbf{4}, 1984, pages 261–281.

\bibitem{HuLau1993}  T.-Y. Hu, K.-S. Lau, Fractal Dimensions and Singularities of the Weierstrass Type Functions,
Transactions of the American Mathematical Society, 1993, \textbf{335}(2), pages 649-665.


\bibitem{WShen2015}  Weixiao Shen, Hausdorff dimension of the graphs of the classical Weierstrass functions, arXiv:1505.03986.



\bibitem{Keller2017} G. Keller, A simpler proof for the dimension of the graph of the classical {W}eierstrass function, Ann. Inst. Poincar\'e, 53(1), 2017, pages 169-181.


\bibitem{KigamiLapidus1993} J. Kigami and M. Lapidus,  {Weyl's problem for the spectral distribution of Laplacians on P.C.F. self-similar fractalss},
{Communications in Mathematical Physics}, \textbf{204}, 2003,
pages 399-444.


\bibitem{Weierstrass} K. Weierstrass, \"{U}ber continuirliche Funktionen eines reellen Arguments, die f\"{u}r keinen Werth des letzteren einen bestimmten Differentialquotienten besitzen, 1967, in Karl Weiertrass Mathematische Werke, Abhandlungen II, Johnson, Gelesen in der Königl. Akademie der Wissenchaften am 18 Juli 1872, \textbf{2}, pages 71-74.

\bibitem{Titschmarsh1939} E. C. Titschmarsh, The theory of functions, Second edition, Oxford University Press, 1939, pages 351-353.


 \bibitem{Hardy} G. H. Hardy, Theorems Connected with Maclaurin's Test for the Convergence of Series, Proc. London Math. Soc., 1911, s2-9 (1), pages 126-144.



 \bibitem{Besicovitch} A. S. Besicovitch, H. D. Ursell, Sets of Fractional Dimensions,  Journal of the London Mathematical Society, \textbf{12} (1), 1937, pages 18-25.

\bibitem{Mandelbrot} B. B. Mandelbrot, Fractals: form, chance, and dimension, San Francisco: Freeman, 1977.


\bibitem{Falconer}  K. Falconer, The Geometry of Fractal Sets, 1985, Cambridge University Press, pages 114-149.


\bibitem{Hunt} B. Hunt, The Hausdorff dimension of graphs of Weierstrass functions, Proc. Amer. Math. Soc., \textbf{126} (3), 1998, pages 791-800.


\bibitem{Baransky} K. Bara\'{n}sky, B. B\'{a}r\'{a}ny, J. Romanowska, On the dimension of the graph of the classical Weierstrass function, Advances in Math., \textbf{265}, 2014, pages 32-59.

\bibitem{Keller} G. Keller, A simpler proof for the dimension of the graph of the classical Weierstrass function, Ann. Inst. Poincar\'e, \textbf{53}(1), 2017, pages 169-181.

\bibitem{Barnsley} M. F. Barnsley, S. Demko, Iterated Function Systems and the Global Construction of Fractals, The Proceedings of the Royal Society of London, \textbf{A}(399), 1985, pages 243-275.

\bibitem{Berry}  M.V. Berry, and Z. V. Lewis,  On the Weierstrass-Mandelbrot function,
Proc. R. Soc. Lond.,  \textbf{A}(370), 1980, pages 459-484.


\bibitem{Sabot1937} C.~Sabot, Existence and uniqueness of diffusions on finitely ramified self-similar fractals, Annales scientifiques de l'\'{E}.N.S. 4 e s\'erie, \textbf{30}(4), 1997, pages 605-673.

\bibitem{Beurling} A. Beurling, J.~Deny, Espaces de Dirichlet. I. Le cas \'el\'ementaire, Acta Mathematica, \textbf{99} (1), 1985, pages 203-224.

\bibitem{Fukushima1994} M. Fukushima, Y. Oshima, and M. Takeda, Dirichlet forms and symmetric Markov processes, 1994, Walter de Gruyter \& Co.

\bibitem{Kigami} J. Kigami, Harmonic Analysis for Resistance Forms, Journal of Functional Analysis, \textbf{204}, 2003, pages 399-444.


\bibitem{Riane}  N. Riane,  2016, Autour du Laplacien sur des domaines pr\'esentant un caract\`ere fractal, M\'emoire de recherche, M2 Math\'ematiques de la mod\'elisation, Universit\'e Pierre et Marie Curie-Paris 6.

\bibitem{DavidRiane} Cl. David et N. Riane, Formes de Dirichlet et fonctions harmoniques sur le graphe de la fonction de Weierstrass, preprint, HAL.




\bibitem{Strichartz2003} R. S. Strichartz, Function spaces on fractals, Journal of Functional Analysis, 198(1), 2003, pages 43-83.

\bibitem{StrichartzLivre2006} R. S. Strichartz, Differential Equations on Fractals, A tutorial, Princeton University Press, 2006.


\bibitem{Hutchinson1971} J. E. Hutchinson, Fractals and self similarity, Indiana University Mathematics Journal 30, 1981, pages 713-747.


\bibitem{StrichartzMeasure} R. S. Strichartz, A. Taylor and T. Zhang, Densities of Self-Similar Measures on the Line, Experimental Mathematics, \textbf{4}(2), 1995, pages 101-128.



\bibitem{Fukushima1994} M. Fukushima and T. Shima, On a spectral analysis for the Sierpinski gasket, Potential Anal., \textbf{1}, 1992, pages 1-3.


\end{thebibliography}
\end{document}